\documentclass{birkjour}

\usepackage{amsmath}
\usepackage{amssymb}
\usepackage{amscd}
\usepackage{latexsym}
\usepackage{epsfig,comment}

\theoremstyle{plain}
\newtheorem{theorem}{Theorem}[section]
\newtheorem{lemma}[theorem]{Lemma}
\newtheorem{proposition}[theorem]{Proposition}
\newtheorem{corollary}[theorem]{Corollary}

\numberwithin{equation}{section}

\theoremstyle{definition}
\newtheorem{definition}[theorem]{Definition}

\newtheorem{remark}[theorem]{Remark}

\newcommand{\rank}{\operatorname{rank}}

\newcommand{\im}{\operatorname{im}}
\newcommand{\id}{\operatorname{id}}

\newcommand{\Hom}{\operatorname{Hom}}
\newcommand{\Ext}{\operatorname{Ext}}

\newcommand{\vac}[1]{ \left| #1 \right> }
\def\Cl{\mathrm{Cl}}

\newcommand{\C}{{\mathbb{C}}}
\newcommand{\Z}{{\mathbb{Z}}}
\newcommand{\N}{{\mathbb{N}}}

\renewcommand{\H}{{\mathcal{H}}}

\renewcommand{\iff}{\Leftrightarrow}

\newcommand{\M}{\mathcal{M}}

\newcommand{\tnv}{\mathrm{tnv}}
\newcommand{\bc}{\mathbf{c}}
\newcommand{\bI}{\mathbf{I}}
\newcommand{\bJ}{\mathbf{J}}
\newcommand{\blambda}{{\boldsymbol\lambda}}

\DeclareMathOperator{\width}{width}
\DeclareMathOperator{\End}{End}

\def\pt{\mathrm{pt}}

\def\K{\mathcal{K}}

\newcommand{\Span}{\mathop{\mathrm{Span}}}

\newcommand{\ag}{\widehat{\mathfrak{g}}}

\begin{document}

\author{Anthony Licata}
\address{%
Department of Mathematics\\
Stanford University\\
Palo Alto, CA\\
USA}
\email{amlicata@math.stanford.edu}

\author{Alistair Savage}
\address{%
Department of Mathematics and Statistics\\
University of Ottawa\\
Ottawa, ON\\
Canada}%
\email{alistair.savage@uottawa.ca}

\thanks{The research of the second author was supported by
an NSERC Discovery Grant.}

\subjclass{Primary: 14C05,17B69}

\title[Vertex operators and geometry]{Vertex operators and the geometry
of\\ moduli spaces of framed torsion-free sheaves}

\begin{abstract}
We define complexes of vector bundles on products of moduli spaces
of framed rank $r$ torsion-free sheaves on $\mathbb{P}^2$.  The top
non-vanishing equivariant Chern classes of the cohomology of these
complexes yield actions of the $r$-colored Heisenberg and Clifford
algebras on the equivariant cohomology of the moduli spaces.  In
this way we obtain a geometric realization of the boson-fermion
correspondence and related vertex operators.
\end{abstract}

\maketitle

\thispagestyle{empty}



\section*{Introduction}

Vertex operators play a vital role in the representation theory of
affine Lie algebras.  The central idea is that instead of
considering individual generators of certain algebras, one should
construct formal generating series of such operators.  In this
formalism, the necessary relations in various constructions become
tractable (for an introduction to these topics, see
\cite{FLM88,LL04}). One of the fundamental building blocks of the
vertex operator calculus is the boson-fermion correspondence (see
\cite{Fre81}), which also plays an important role in mathematical
physics.

The relationship between Hilbert schemes of points on surfaces and
the representation theory of Heisenberg algebras has been well
studied.  In particular, Nakajima \cite{Nak99} defined an action of
the Heisenberg algebra on the cohomology of such Hilbert schemes,
endowing this cohomology with the structure of bosonic Fock space.
Relations to vertex operators were demonstrated by Grojnowski
\cite{Gro96}, Lehn \cite{Leh99}, Li-Qin-Wang \cite{LQW02}, and
Carlsson-Okounkov \cite{CO08}. The construction of the Heisenberg
algebra was then extended to the case of equivariant cohomology by
Vasserot \cite{Vas01}.  In this setting, connections to vertex
operator constructions appeared in \cite{LQW04}.  These
relationships suggest that vertex operators themselves should have
geometric descriptions in terms of the Hilbert scheme and related
moduli spaces.  In this paper, we give such a description (see
Section~\ref{sec:BFC-geometric-realization}).

In \cite{Sav06b}, the second author pointed out that the
$\C^*$-fixed points of the Hilbert schemes of points in $\C^2$ can
be naturally identified with quiver varieties of type $A_\infty$. By
a general theory of Nakajima \cite{Nak94,Nak98}, the homology (or
equivariant cohomology) of these quiver varieties can be given the
structure of a level one representation of $\mathfrak{gl}_\infty$.
Thus, the localization theorem, which relates the equivariant
cohomology of a space with that of its fixed point set, yields a
geometric setting for the Fock space realization of the basic
representation of $\mathfrak{gl}_\infty$.  A similar idea was used
by the first author in \cite{Lic07} to construct the action of the
Clifford algebra instead of the Lie algebra $\mathfrak{gl}_\infty$.

While these approaches provide a geometric interpretation of the
boson-fermion correspondence, the action of the Clifford algebra is
via correspondences on a finite fixed point set of a torus action
while the action of the Heisenberg algebra is defined via
correspondences on the entire Hilbert schemes.  It would be
interesting to put the action of the Clifford algebra on the same
footing as that of the Heisenberg algebra.  That is, one would like
to have globally defined operators, rather than just operators
living on torus fixed point sets.  In this way, the boson-fermion
correspondence would appear as the action of both the Clifford and
Heisenberg algebras on the same space -- the equivariant cohomology
of Hilbert schemes or their higher rank generalizations. Apart from
being satisfying in its own right, this should be an important step
in the program of understanding the entire vertex operator calculus
in terms of the moduli space geometry.

In the current paper, we accomplish the aforementioned task as well
as the generalization to higher rank. Namely, we define actions of
the $r$-colored Clifford and Heisenberg algebras on the equivariant
cohomology of moduli spaces of framed rank $r$ torsion-free sheaves
on $\mathbb{P}^2$. The action is defined using equivariant Chern
classes of vector bundles defined via complexes on products of these
moduli spaces, yielding a new and uniform realization of the
Heisenberg and Clifford algebra actions.  This method of defining
operators was inspired by the work of Carlsson and Okounkov
\cite{CO08}, who realized a large class of vertex operators as Chern
classes of virtual vector bundles on Hilbert schemes on an arbitrary
(quasi-)projective surface. Indeed, our vector bundles have a
natural description in the language of sheaves, and in the case of
framed rank one sheaves on $\mathbb{P}^2$, they are closely related
to the virtual vector bundles considered by Carlsson and Okounkov
for the surface ${\C^2}$ (see Section~\ref{sec:sheaf-interpretation}
for details).

Our extension to higher rank is important in that the $r$-colored
Fock spaces play a key role in other constructions in the theory of
vertex operators, such as the vertex operator realization of the
basic representations of affine Lie algebras.  In particular, we
predict that complexes similar to those defined in this paper will
lead to moduli space constructions of all vertex operator
realizations (including the homogeneous and principle realizations)
of the basic representations.  Since our constructions are defined
via complexes of vector bundles, we expect this method to extend
naturally to K-theory and even to a full categorification of
representations of the Heisenberg algebras, Clifford algebras, and
affine Lie algebras. We also anticipate obtaining natural
interpretations of level-rank duality as a certain compatibility
between different geometrically defined operators on the equivariant
cohomology of moduli spaces.

The organization of this paper is as follows.  In
Section~\ref{sec:Fock-spaces} we review the $r$-colored bosonic and
fermionic Fock spaces.  The moduli spaces of framed rank $r$
torsion-free sheaves and various results about their equivariant
cohomology are discussed in Section~\ref{sec:moduli-spaces}. There
we also discuss how one can define operators on this equivariant
cohomology via equivariant Chern classes of virtual vector bundles.
In Section~\ref{sec:geometric-operators} we define the geometric
Clifford and Heisenberg operators, leaving the proofs of the main
theorems to Sections~\ref{sec:proof-clifford}
and~\ref{sec:proof-heisenberg}. We make some concluding remarks as
well as indications of direction of further research in
Section~\ref{sec:future}.

The authors would like to thank I. Frenkel, H. Nakajima, and W. Wang
for helpful comments on an earlier version of this paper. They would
also like to thank H. Nakajima for helping point out the explicit
relationship between some of the constructions appearing in this
paper and those in \cite{CO08} (see
Section~\ref{sec:sheaf-interpretation}) and D. Maulik for bringing
to their attention an omission in an earlier version of this paper.


\section{Fock spaces} \label{sec:Fock-spaces}

In this section we introduce the $r$-colored bosonic and fermionic
Fock spaces.  Further exposition can be found in the references
mentioned in Section~\ref{sec:boson-fermion-corr}.

\subsection{Bosonic Fock space} \label{sec:bosonic-fock-space}

A \emph{partition} is an infinite sequence
\[
    \lambda = (\lambda_1, \lambda_2, \dots)
\]
with a finite number of non-zero terms and such that $\lambda_1 \ge
\lambda_2 \ge \dots$.  We identify a finite sequence $(\lambda_1,
\dots, \lambda_n)$ with the infinite sequence obtained by setting
$\lambda_i=0$ for $i>n$.  Let $\mathcal{P}$ denote the set of all
partitions. For a partition $\lambda = (\lambda_1, \lambda_2,
\dots)$, let $\lambda^t= ((\lambda^t)_1, (\lambda^t)_2, \dots)$
denote the dual (transposed) partition. That is,
\[
    (\lambda^t)_i = \# \{j\ |\ \lambda_j \ge i\}.
\]
The largest integer $i$ such that $\lambda_i \ne 0$ is called the
\emph{length} of $\lambda$ (if $\lambda_i=0$ for all $i$, the length
of $\lambda$ is zero). Denote the \emph{size} of the partition
$\lambda$ by $|\lambda| = \sum_{i} \lambda_i$.   Following
\cite{NY05}, we will identify a partition $\lambda$ with the set
\[
\{(i,j)\ |\ i \ge 1,\ 1 \le j \le \lambda_i\} \subset \N_+^2,
\]
where $\N_+$ denotes the set of (strictly) positive integers. When
we write $s \in \lambda$, we are referring to this identification.
Note that this convention differs from the usual English or French
notations for Young diagrams.

For $s =(i,j) \in \N_+^2$, we define the \emph{arm} and \emph{leg
length of $s$ in $\lambda$} by
\[
    a_{\lambda}(s) = \lambda_i - j \quad \text{and} \quad
    l_{\lambda}(s) = (\lambda^t)_j - i
\]
respectively.  Note that these numbers are negative if $s \not \in
\lambda$.

Given a pair of partitions $\lambda, \mu$ and $s \in \N_+^2$, define
the \emph{relative hook length} of $s$ to be
\[
h_{\lambda,\mu}(s) = a_\lambda(s) + l_\mu(s) + 1
\]
and the \emph{relative hook number}
\[
    h_{\lambda,\mu} = \prod_{s \in \lambda} h_{\lambda,\mu}(s).
\]
The usual \emph{hook length of $s$ in $\lambda$} is then
\[
    h_\lambda(s) = h_{\lambda, \lambda}(s)
\]
and the \emph{hook number $h_\lambda$ of $\lambda$} is
\[
    h_\lambda =  h_{\lambda, \lambda}.
\]

\begin{lemma} \label{lem:hook=0}
Suppose $\lambda \ne \mu$ are partitions with
\[
\lambda_k > \mu_k,\quad \text{and}\quad \lambda_i = \mu_i,\quad 1
\le i < k.
\]
Then $h_{\lambda,\mu}((k,\lambda_k))=0$.
\begin{proof}
We have
\[
h_{\lambda,\mu}((k,\lambda_k)) = a_\lambda((k,\lambda_k)) +
l_\mu((k,\lambda_k)) + 1 = 0 + (-1) + 1 = 0.
\]
\end{proof}
\end{lemma}

Let $\Lambda$ denote the ring of symmetric functions with complex
coefficients. We let $\{p_n\}_{n \in \N}$ and
$\{s_\lambda\}_{\lambda \in \mathcal{P}}$ denote the power sums and
Schur functions respectively.  The Schur functions form a basis of
$\Lambda$ and we also have
\[
    \Lambda = \C[p_1, p_2, \dots].
\]
We fix the bilinear form $\langle \cdot, \cdot \rangle$ on $\Lambda$
with respect to which the basis of Schur functions is orthonormal.

Fix $r \in \N_+$ and define the \emph{$r$-colored oscillator
algebra} to be the Lie algebra
\[
    \mathfrak{s} = \bigoplus_{m \in \Z,\, l \in \{1, \dots,
    r\}} \C p^l(m) \oplus \C K
\]
with commutation relations
\[
    [\mathfrak{s}, K] = [\mathfrak{s},p^k(0)] = 0,\quad [p^k(m), p^l(n)] =
    \frac{1}{m} \delta_{m,-n} \delta_{k,l} K,\ m \ne 0.
\]
The subalgebra spanned by $p^l(n)$, $l \in \{1, \dots, r\}$, $n \in
\Z \setminus \{0\}$, and $K$ is an $r$-colored infinite-dimensional
Heisenberg algebra.

The oscillator algebra has a natural representation on the
\emph{$r$-colored bosonic Fock space}
\[
    \mathbf{B} = B^{\otimes r}, \quad B \stackrel{\text{def}}{=}
    \C[p_1,p_2,\dots;q,q^{-1}] \cong \Lambda \otimes_\C \C[q,q^{-1}]
\]
given by:
\begin{align*}
    p^l(m) &\mapsto \id^{\otimes (l-1)} \otimes \frac{\partial}{\partial
    p_m} \otimes \id^{\otimes (r-l)},\ m > 0,\\
    p^l(-m) &\mapsto \id^{\otimes (l-1)}
    \otimes \frac{1}{m} p_m \otimes \id^{\otimes (r-l)},\ m>0, \\
    p^l(0) &\mapsto \id^{\otimes (l-1)} \otimes q \frac{\partial}{\partial q}
    \otimes \id^{\otimes (r-l)},\quad K \mapsto \id.
\end{align*}
Note that
\begin{align*}
    \mathbf{B} &= \bigoplus_{\bc \in \Z^r, n \in \N}
    \mathbf{B}^\bc_n,\\
    \mathbf{B}^\bc_n &\stackrel{\text{def}}{=}
    \left\{ \left. (q^{c^1} f_1, \dots, q^{c^r} f_r)\ \right|\ f_\alpha \in \Lambda
    \, \forall \alpha,\, \sum_{\alpha = 1}^ r \deg
    f_\alpha = n\right\}.
\end{align*}
We fix the bilinear form on $\mathbf{B}$ given by
\[
    \left< (q^{c^1} f_1, \dots ,q^{c^r} f_r), (q^{d^1} g_1, \dots, q^{d^r}
    g_r) \right> = \delta_{c^1, d^1} \cdots \delta_{c^r, d^r} \left< f_1,
    g_1 \right> \cdots \left< f_r, g_r \right>
\]
for $f_\alpha, g_\alpha \in \Lambda$, $\alpha \in \{1, \dots, r\}$.
It is easily verified that the operators $p^l(m)$ and $p^l(-m)$ are
adjoint for $m \in \Z$.

For an $r$-tuple $\blambda = (\lambda^1, \dots, \lambda^r)$ of
partitions, define $h_{\blambda} = \prod_{\alpha=1}^r
h_{\lambda^\alpha}$ and $|\blambda| = \sum_{\alpha=1}^r
|\lambda^\alpha|$.  We call $r$ the \emph{rank} of the oscillator
algebra and bosonic Fock space.


\subsection{Fermionic Fock space} \label{sec:fermionic-fock-space}
An infinite expression of the form
\[
    i_1 \wedge i_2 \wedge i_3 \wedge \dots,
\]
where $i_1, i_2, \dots$ are integers satisfying
\[
    i_1 > i_2 > i_3 > \dots,\ i_n = i_{n-1} - 1 \text{ for } n \gg 0,
\]
is called a \emph{semi-infinite monomial}.  Let $F$ be the complex
vector space with basis consisting of all semi-infinite monomials.
Then $F$ is called \emph{fermionic Fock space}. Let
\[
    \vac{c} = c \wedge (c-1) \wedge (c-2) \wedge \dots
\]
be the \emph{vacuum vector of charge $c$}.  We say that a
semi-infinite monomial has \emph{charge} $c$ if it differs from
$\vac{c}$ at only a finite number of places.  Thus $I = i_1 \wedge
i_2 \wedge \dots$ is of charge $c$ if $i_k = c - k + 1$ for $k \gg
0$.  Let $F^{c}$ denote the linear span of all semi-infinite
monomials of charge $c$.  Then we have the charge decomposition
\[
    F = \bigoplus_{c \in \Z} F^{c}.
\]

A semi-infinite monomial $I = i_1 \wedge i_2 \wedge \dots$
determines a partition $\lambda(I) \in \mathcal{P}$ by
\[
    i_k = (c(I)-k+1) + \lambda(I)_k,\ k \in \N_+.
\]
This gives a bijection $I \mapsto (\lambda(I), c(I))$ between the
set of all semi-infinite monomials and the set $\mathcal{P} \times
\Z$. We define the \emph{energy} of $I$ to be $|I|=|\lambda(I)|$.
Let $F_j^{c}$ denote the linear span of all semi-infinite monomials
of charge $c$ and energy $j$. We then have the energy decomposition
\[
F^{c} = \bigoplus_{j \in \N} F_j^{c}.
\]

For $j \in \Z$, define the \emph{wedging} and \emph{contracting}
operators $\psi(j)$ and $\psi(j)^*$ on $F$ by:
\begin{align*}
    \psi(j)(i_1 \wedge i_2 \wedge \dots) &=
    \begin{cases}
    0 & \text{if $j = i_s$ for some $s$,} \\
    (-1)^s i_1 \wedge \dots \wedge i_s
    \wedge j \wedge i_{s+1} \wedge \dots &
    \text{if $i_s > j > i_{s+1}$}.
    \end{cases} \\
    \psi(j)^*(i_1 \wedge i_2 \wedge \dots) &=
    \begin{cases}
    0 & \text{if $j \ne i_s$ for all $s$,} \\
    (-1)^{s-1} i_1 \wedge \dots \wedge
    i_{s-1} \wedge i_{s+1} \wedge \dots &
    \text{if $j=i_s$}.
    \end{cases}
\end{align*}

Fix $r \in \N_+$.  We define the \emph{$r$-colored Clifford algebra}
$\Cl$ to be the algebra with generators $\psi^l(j), \psi^l(j)^*$ $j
\in \Z,\, l \in \{1,\dots,r\}$ and relations
\begin{gather*}
    \{\psi^l(j), \psi^l(i)^*\} = \delta_{ji},\quad
    \{\psi^l(j), \psi^l(i)\} = 0 =
    \{\psi^l(j)^*, \psi^l(i)^*\}, \\
    [\psi^l(j), \psi^k(i)] = [\psi^l(j), \psi^k(i)^*] =
    [\psi^l(j)^*, \psi^k(i)^*] = 0,\quad l \ne k,
\end{gather*}
where $\{a,b\} = ab + ba$. Note that we choose the convention that
generators of different colors commute since it turns out that this
is more natural from the geometric viewpoint. If one prefers to
consider a Clifford algebra where generators of different colors
anti-commute, it is only necessary to insert an appropriate sign in
the geometric operators to be defined. We define \emph{$r$-colored
fermionic Fock space} to be the vector space
\[
    \mathbf{F} = F^{\otimes r}.
\]
We call an $r$-tuple $\bI = (I^1, \dots, I^r)$ of semi-infinite
monomials an \emph{$r$-colored semi-infinite monomial}.  These form
a basis of $\mathbf{F}$ and we let $\langle \cdot, \cdot \rangle$
denote the bilinear form on $\mathbf{F}$ for which this basis is
orthonormal. For $l \in \{1, 2, \dots, r\}$ and $j \in \Z$, the maps
\begin{align*}
    \psi^l(j) &= \id^{\otimes (l-1)} \otimes \psi(j) \otimes
    \id^{\otimes (r-l)},\quad \text{and} \\
    \psi^l(j)^* &= \id^{\otimes (l-1)} \otimes \psi(j)^* \otimes
    \id^{\otimes (r-l)}
\end{align*}
define a representation of $\Cl$ on $\mathbf{F}$.  Note that
\[
    \psi^l(j) (F^\bc) \subseteq F^{\bc+1_l},\quad \psi^l(j)^*(F^\bc)
    \subseteq F^{\bc-1_l},
\]
where $1_l \in \Z^r$ has a one in the $l$th position and a zero
everywhere else.  The operators $\psi^l(j)$ and $\psi^l(j)^*$ are
called \emph{free fermions}. One can check directly that $\psi^l(j)$
and $\psi^l(j)^*$ are adjoint with respect to the bilinear form
$\langle \cdot, \cdot \rangle$ and that $\mathbf{F}$ is an
irreducible $\Cl$-module.  We call $r$ the \emph{rank} of the
Clifford algebra and fermionic Fock space.

We have the charge decomposition
\[
    \mathbf{F} = \bigoplus_{\bc \in \Z^r} \mathbf{F}^{\bc},\quad
    \mathbf{F}^{\bc} = F^{c^1} \otimes \dots \otimes F^{c^r},
\]
and energy decomposition
\[
    \mathbf{F}^{\bc} = \bigoplus_{n \in \N}
    \mathbf{F}^{\bc}_n,\quad \mathbf{F}^{\bc}_n = \bigoplus_{j^1 +
    \dots + j^r = n}
    F^{c^1}_{j^1} \otimes \dots \otimes F^{c^r}_{j^r}.
\]

In order to simplify notation in what follows, for a semi-infinite
monomial $I$, we define
\[
    a_I(s) = a_{\lambda(I)}(s),\quad l_I(s) = l_{\lambda(I)}(s).
\]
We define $h_I(s)$, $h_I$, $h_{I,J}(s)$ and $h_{I,J}$ similarly. For
an $r$-tuple $\bI = (I^1, \dots, I^r)$ of semi-infinite monomials,
we define $\blambda(\bI) = (\lambda(I^1), \dots, \lambda(I^r))$,
$h_\bI = h_{\blambda(\bI)}$, $|\bI| = |\blambda(\bI)|$, and
$\bc(\bI) = (c(I^1), \dots, c(I^r))$.


\subsection{The boson-fermion correspondence}
\label{sec:boson-fermion-corr}

The boson-fermion correspondence is a precise relationship between
bosonic and fermionic Fock space (see \cite{Fre81}).  It uses vertex
operators to express bosons in terms of fermions (bosonization) and
fermions in terms of bosons (fermionization).  We do not present the
details here but instead refer the reader to the expository
presentations of this topic found in \cite[Chapter~14]{K} and the
introduction to \cite{tKvdL91}.  In
Section~\ref{sec:BFC-geometric-realization} we present a geometric
version of the boson-fermion correspondence.


\section{Moduli spaces of framed torsion-free sheaves on $\mathbb{P}^2$}
\label{sec:moduli-spaces}

In this section we discuss our main object of study: the moduli
space of framed torsion-free sheaves on the complex projective plane
$\mathbb{P}^2$ and its equivariant cohomology.  For background on
equivariant cohomology, especially equivariant Chern classes, we
refer the reader to \cite[Chapter 9]{CK99}.

\subsection{The moduli space $\M(r,n)$} \label{sec:Mrn}

Let $\M(r,n)$ be the moduli space of framed torsion-free sheaves on
$\mathbb{P}^2$ with rank $r$ and second Chern class $c_2=n$. More
precisely, points of $\M(r,n)$ are isomorphism classes of pairs
$(E,\Phi)$ where $E$ is a torsion free sheaf with $\rank E = r$ and
$c_2(E)=n$ which is locally free in a neighborhood of $l_\infty$,
and $\Phi : E_{l_\infty} \xrightarrow{\cong}
\mathcal{O}_{l_\infty}^r$ is a framing at infinity. Here $l_\infty =
\{[0 : z_1 : z_2] \in \mathbb{P}^2\}$ is the line at infinity. The
existence of a framing implies $c_1(E)=0$.

There is another description of $\M(r,n)$ (basically due to Barth
\cite{Bar77}) which will be important for us.  Fix $r \in \N_+$, $n
\in \N$ and let $V = \C^n$, $W = \C^r$.
\begin{proposition}[{\cite[Theorem~2.1]{Nak99}}]
There exists an isomorphism of algebraic varieties
\[
\M(r,n) \cong \{(A,B,i,j)\ |\ [A,B]+ij=0,\, (A,B,i,j) \text{ is
stable}\}/GL(V)
\]
where $A,B \in \End V$, $i \in \Hom (W,V)$, $j \in \Hom(V,W)$ and
the $GL(V)$-action is given by
\[
g (A,B,i,j) = (g A g^{-1}, g B g^{-1}, gi, jg^{-1}),\quad g \in
GL(V).
\]
We say $(A,B,i,j)$ is \emph{stable} if there exists no proper
subspace $S \subsetneq V$ such that $A(S) \subseteq S$, $B(S)
\subseteq S$ and $\im i \subseteq S$.
\end{proposition}
We call the quadruple $(A,B,i,j)$ in the above proposition the
\emph{ADHM data} corresponding to the framed sheaf $(E, \Phi)$. We
note that $\M(1,n)$ is isomorphic to the Hilbert scheme of $n$
points in $\C^2$.  Thus, the moduli space $\M(r,n)$ can be thought
of as a higher rank version of this Hilbert scheme.


\subsection{Torus actions on $\M(r,n)$} \label{sec:torus-action}

Fix the torus $T = (\C^*)^r \times \C^*$.  Elements of $(\C^*)^r$
will be written
\[
    e = (e_1,\hdots,e_r).
\]
We denote the one-dimensional $T$-modules
\[
    (e,t) \mapsto e_i \quad \text{and} \quad (e,t)\mapsto t
\]
by $e_i$ and $t$ respectively and the tensor product of such modules
by juxtaposition. If $\pt$ is the space consisting of a single point
with the trivial $T$-action, then $H_T^*(\pt) = \C[b_1, b_2, \dots,
b_r, \epsilon]$ where $\epsilon$ and $b_i$, $i \in \{1, 2, \dots,
r\}$, are the second Chern classes of the one-dimensional
$T$-modules $t$ and $e_i$ respectively.  They are elements of degree
2. We fix a decomposition $W = \bigoplus_{\alpha=1}^r W^\alpha$
where $W^\alpha \cong \C$ for all $\alpha$.  Then $(\C^*)^r$ acts on
$W$ by
\[
    e \mapsto e_1 \id_{W^1} \oplus \dots \oplus e_r \id_{W^r}.
\]
For $\bc \in \Z^r$, let $\M_{\bc}(r,n)$ denote the moduli space with
the torus action
\[
    (e,t) \star_{\bc}(A,B,i,j) = (tA, t^{-1}B, ie^{-1} t^{-\bc},
    e t^{\bc} j),
\]
where
\[
    t^\bc = (t^{c^1}, \dots, t^{c^r}) \in (\C^*)^r.
\]
Note that the underlying variety of $\M_\bc(r,n)$ is independent of
$\bc$.  It is only the $T$-action that changes.

The $T$-fixed points $\M_\bc(r,n)^T$ are in natural bijection with
$r$-tuples of semi-infinite monomials $\bI = (I^1,\hdots,I^r)$, with
$\bc(\bI) = \bc$ and $|\blambda(\bI) | = n$ (see
\cite[Proposition~2.9]{NY05}). We shall identify $T$-fixed points
with such $r$-tuples in what follows.  Let $\bI = (I^1, \dots, I^r)
\in \M_\bc(r,n)^T$ and let $\mathcal{T}_\bI$ denote the tangent
space to $\M_\bc(r,n)$ at the point $\bI$. Then $T$ acts on
$\mathcal{T}_\bI$, which decomposes into one-dimensional
$T$-modules.

\begin{proposition}
As a $T$-module, $\mathcal{T}_\bI$ is given by
\[
    \mathcal{T}_\bI = \bigoplus_{\alpha,\beta = 1}^r \left(
    e_\beta e_\alpha^{-1} t^{c(I^\beta)-c(I^\alpha)}
    \left( \bigoplus_{s \in \lambda(I^\alpha)} t^{-h_{I^\alpha,
    I^\beta}(s)} \oplus \bigoplus_{s \in \lambda(I^\beta)}
    t^{h_{I^\beta, I^\alpha}(s)} \right) \right).
\]
\end{proposition}
\begin{proof}
This follows from \cite[Theorem~2.11]{NY05} after replacing $t_1$ by
$t$, $t_2$ by $t^{-1}$, and $e_{\alpha}$ by $e_\alpha
t^{c(I^\alpha)}$ everywhere.
\end{proof}

It will be convenient to consider the splitting $\mathcal{T}_\bI =
\mathcal{T}_\bI^- \oplus \mathcal{T}_\bI^+$ where
\begin{align*}
    \mathcal{T}_\bI^- &= \bigoplus_{\alpha, \beta =1}^r \left( e_\beta
    e_\alpha^{-1} t^{c(I^\beta) - c(I^\alpha)} \bigoplus_{s \in
    \lambda(I^\alpha)} t^{-h_{I^\alpha, I^\beta}(s)} \right), \\
    \mathcal{T}_\bI^+ &= \bigoplus_{\alpha, \beta =1}^r \left(
    e_\alpha
    e_\beta^{-1} t^{c(I^\alpha) - c(I^\beta)} \bigoplus_{s \in
    \lambda(I^\alpha)} t^{h_{I^\alpha, I^\beta}(s)} \right).
\end{align*}

\begin{lemma} \label{lem:Euler-classes-tangent-spaces}
The equivariant Euler classes of $\mathcal{T}_\bI^-$ and
$\mathcal{T}_\bI^+$ in the $T$-equivariant cohomology of a point are
given by
\begin{align*}
    e_T(\mathcal{T}_\bI^-) &=
    \prod_{\alpha, \beta = 1}^r \prod_{s \in \lambda(I^\alpha)}
    (b_\beta - b_\alpha + (c(I^\beta) - c(I^\alpha) - h_{I^\alpha,
    I^\beta}(s))\epsilon),\\
    e_T(\mathcal{T}_\bI^+) &= \prod_{\alpha, \beta=1}^r
    \prod_{s\in \lambda(I^\alpha)}
    (b_\alpha - b_\beta + (c(I^\alpha) - c(I^\beta) + h_{I^\alpha,
    I^\beta}(s))\epsilon) \\
    &= (-1)^{r|\bI|}e_T(\mathcal{T}_\bI^-).
\end{align*}
\end{lemma}
\begin{proof}
This follows easily from the definition of $\mathcal{T}_\bI^-$ and
$\mathcal{T}_\bI^+$.
\end{proof}


\subsection{Bilinear form and the space $\mathbf{A}$} \label{sec:bilinear-form}

Let
\[
\mathcal{H}_T^*(\M_\bc(r,n)) = H_T^*(\M_\bc(r,n)) \otimes_{\C[b_1,
\dots, b_r, \epsilon]} \C(b_1, \dots, b_r, \epsilon)
\]
denote the localized equivariant cohomology.  Let
\[
    i: \M_\bc(r,n)^T \hookrightarrow \M_\bc(r,n)
\]
denote inclusion, and let
\[
    p:  \M_\bc(r,n)^T \twoheadrightarrow \{\pt\}
\]
be the projection to a point.  The real dimension of $\M_\bc(r,n)$
is $4rn$.  We define a bilinear form $ \langle \cdot, \cdot
\rangle_{n, \bc} $ on the middle degree localized equivariant
cohomology $\mathcal{H}^{2rn}_T(\M_\bc(r,n))$ by
\[
    \langle a,b\rangle_{n, \bc} = (-1)^{rn}p_*(i_*)^{-1}(a \cup b),
\]
where $i_*$ is invertible by the localization theorem.  We then
extend this to a bilinear form
\[
    \langle \cdot, \cdot \rangle = \bigoplus_{n,\bc} \langle \cdot, \cdot
    \rangle_{n, \bc}
\]
on $\bigoplus_{n,\bc} \mathcal{H}^{2rn}_T(\M_\bc(r,n))$ by declaring
classes from different summands to be orthogonal.  For $\bI \in
\M_\bc(r,n)^T$, let
\[
    [\bI] = \frac{i_*(1_{\bI})}{e_T(\mathcal{T}_{\bI}^-)} \in
    \mathcal{H}_T^{2r|\bI|}(\M_{\bc(\bI)}(r,|\bI|)).
\]
Here $1_\bI$ denotes the unit in the $T$-equivariant cohomology ring
of a point and the denominator $e_T (\mathcal{T}_\bI^-)$ is to be
interpreted as an element in this ring.  By the localization theorem
in equivariant cohomology, the elements $[\bI]$ form a $\C(b_1, b_2,
\dots, b_r, \epsilon)$-basis of $\bigoplus_{n,\bc}
\mathcal{H}^*_T(\M_\bc(r,n))$.

\begin{proposition}
The classes $\{[\bI]\}$ are orthonormal.
\end{proposition}
\begin{proof}
For $\bI,\bJ \in \M_{\bc}(r,n)^T$ we compute
\[
    \langle [\bI],[\bJ] \rangle = (-1)^{rn}p_*(i_*)^{-1}([\bI]
    \cup [\bJ]) =
    (-1)^{rn} p_* (i_*)^{-1} \left( \frac{i_*(1_{\bI})}{e_T(\mathcal{T}_{\bI}^-)}
    \cup \frac{i_*(1_{\bJ})}{e_T(\mathcal{T}_{\bJ}^-)} \right)
\]
which is clearly $0$ unless $\bI = \bJ$.  If $\bI = \bJ$, then by
the projection formula the above is equal to
\[
    (-1)^{rn}p_* \left(
    \frac{1_{\bI}}{e_T(\mathcal{T}_{\bI}^-)} \cup \left(
    \frac{i^*i_*(1_{\bI})}{e_T(\mathcal{T}_{\bI}^-)}\right) \right) =
    (-1)^{rn} \frac{e_T(\mathcal{T}_{\bI})}{e_T(\mathcal{T}_{\bI}^-)
    e_T(\mathcal{T}_{\bI}^-)}=
1.
\]
Note that the last equality in the above line follows from our
splitting of $\mathcal{T}_\bI$, which gave
\[
    e_T(\mathcal{T}_{\bI}) = e_T(\mathcal{T}_{\bI}^-)e_T(\mathcal{T}_{\bI}^+) =
    (-1)^{rn}e_T(\mathcal{T}_{\bI}^-)e_T(\mathcal{T}_{\bI}^-).
\]
\end{proof}

A priori, the classes $[\bI]$ are elements of the localized
equivariant cohomology $\mathcal{H}^{2rn}_T(\M_\bc(r,n))$. However,
despite the division by the equivariant Euler class in the
definition, we expect that they lie in the nonlocalized equivariant
cohomology.

Let
\begin{align*}
    A_\bc(r,n) &= \mathrm{Span}_{\C} \{[\bI]\ |\ \bI \in \M_\bc(r,n)^T\} \subset
    \mathcal{H}^{2rn}_T(\M_{\bc}(r,n)), \\
    \mathbf{A} &= \bigoplus_{\bc \in \Z^r,\, n \in \N}
    A_\bc(r,n).
\end{align*}
Recall that the set $\{[\bI]\ |\ \bI \in \M_\bc(r,n)^T\}$ is a
$\C(b_1, \dots, b_r, \epsilon)$-basis of
\[ \textstyle
  \bigoplus_{n, \bc} \mathcal{H}_T^*(\M_\bc(r,n)).
\]
Thus $\mathbf{A}$ is a full $\C$-lattice in this space.

\begin{corollary} \label{cor:bilinear-form-on-A}
The restriction of $\langle \cdot, \cdot \rangle$ to $\mathbf{A}$ is
non-degenerate and $\C$-valued.
\end{corollary}


\subsection{Operators on equivariant cohomology}
\label{sec:operators} We define a bilinear form on the
localized equivariant cohomology of a product of moduli spaces
$\M_{\bc}(r,n_1)\times \M_{\mathbf{d}}(r,n_2)$ in a matter analogous
to that introduced in Section~\ref{sec:bilinear-form}.  Namely,
\[
    \langle a,b \rangle_{n_1, \bc, n_2, \mathbf{d}} =
    (-1)^{rn_2}p_*((i_1 \times i_2)_*)^{-1}(a \cup b),
\]
where $i_1$ and $i_2$ are the inclusions of the $T$-fixed points
into the first and second factors respectively.  We extend this to a
bilinear form
\[
    \langle \cdot, \cdot \rangle = \bigoplus_{n_1, \bc, n_2,
    \mathbf{d}} \langle \cdot, \cdot \rangle_{n_1, \bc, n_2,
    \mathbf{d}}
\]
on $\bigoplus_{n_1, \bc, n_2, \mathbf{d}} \H_T^{2r(n_1+n_2)}
(\M_{\bc}(r,n_1)\times \M_{\mathbf{d}}(r,n_2))$ by declaring
elements in different summands to be orthogonal.

If $\alpha \in \mathcal{H}^{2r(n_1 + n_2)}_T(\M_\bc(r,n_1) \times
\M_{\mathbf{d}}(r,n_2))$, then $\alpha$ defines an operator
\[
    \alpha: \mathcal{H}_T^{2rn_1}(\M_{\bc}(r,n_1))
    \longrightarrow \mathcal{H}^{2rn_2}_T(\M_{\mathbf{d}}(r,n_2))
\]
by using the bilinear form to define structure constants:
\[
    \langle \alpha x,y \rangle_{n_2, \mathbf{d}} \stackrel{\text{def}}{=}
    \langle x\otimes y, \alpha \rangle_{n_1, \bc, n_2, \mathbf{d}}.
\]
In particular, if $E$ is a $T$-equivariant vector bundle on
$\M_{\bc}(r,n_1) \times \M_{\mathbf{d}}(r,n_2)$ and $\beta \in
\mathcal{H}^{2l}_T(\M_\bc(r,n_1) \times \M_{\mathbf{d}}(r,n_2))$
then $\beta \cup c_{r(n_1+n_2)-l}(E)$ defines such an operator.  We
note that in the sequel, our choices will yield operators that
restrict to non-localized equivariant cohomology.

The following lemma is the basic localization tool we will use to
compute the action of our geometric Heisenberg and Clifford
operators.

\begin{lemma} \label{lem:structure-constants}
Suppose $\bI \in \M_{\bc}(r,n_1)^T$ and $\bJ \in
\M_{\mathbf{d}}(r,n_2)^T$.  Let $E$ be a $T$-equivariant vector
bundle on $\M_{\bc}(r,n_1) \times \M_{\mathbf{d}}(r,n_2)$ and $\beta
\in \mathcal{H}^{2l}_T(\M_\bc(r,n_1) \times \M_{\mathbf{d}}(r,n_2))$
for some $l \in \N$. Then
\[
    \langle \beta \cup c_{r(n_1+n_2)-l}(E) [\bI] , [\bJ] \rangle =
    \frac{\beta_{\bI, \bJ} \cup c_{r(n_1+n_2)-l}(E_{(\bI,\bJ)})}
    {e_T(\mathcal{T}_\bI^-) e_T(\mathcal{T}_\bJ^+)},
\]
where $c_{r(n_1+n_2)-l}(E_{(\bI,\bJ)})\in H^*_T(\pt) =
\C[b_1,\hdots,b_r,\epsilon]$ is the polynomial given by the
equivariant Chern class of the fiber of $E$ over the point
$(\bI,\bJ)$ and $\beta_{\bI, \bJ} = i_{\bI,\bJ}^* (\beta)$ where
$i_{\bI,\bJ} : (\bI, \bJ) \hookrightarrow \M_\bc(r,n_1) \times
\M_{\mathbf{d}}(r,n_2)$ is the inclusion of the fixed point.

\begin{proof}
\begin{align*}
    \langle \beta &\cup c_{r(n_1+n_2)-l}(E) [\bI] , [\bJ] \rangle \\
    &= (-1)^{rn_2} p_* ((i_1 \times i_2)_*)^{-1}
    \left([\bI] \otimes [\bJ]  \cup \beta \cup c_{r(n_1+n_2)-l}(E)\right) \\
    &=(-1)^{rn_2} p_* ((i_1 \times i_2)_*)^{-1}
    \left((i_1 \times i_2)_*(1_{(\bI,\bJ)})  \cup
    \frac{\beta \cup c_{r(n_1+n_2)-l}(E)}{e_T(\mathcal{T}_{\bI}^-)
    e_T(\mathcal{T}_{\bJ}^-)}\right).
\end{align*}
where the last line just used the definition of $[\bI]$ and $[\bJ]$.
By the projection formula and the functoriality of Chern classes,
this is equal to
\begin{align*}
    (-1)^{rn_2} p_* &\left(1_{(\bI,\bJ)} \cup \frac{(i_1\times i_2)^*(\beta \cup
    c_{r(n_1+n_2)-l}(E))}
    {e_T(\mathcal{T}_{\bI}^-)e_T(\mathcal{T}_{\bJ}^-)}\right) \\
    &= (-1)^{rn_2} p_* \left(1_{(\bI,\bJ)}  \cup \frac{(i_1\times
    i_2)^*(\beta) \cup
    c_{r(n_1+n_2)-l}((i_1 \times i_2)^*E))}
    {e_T(\mathcal{T}_{\bI}^-)e_T(\mathcal{T}_{\bJ}^-)}\right) \\
    &= \frac{\beta_{\bI, \bJ} \cup c_{r(n_1+n_2)-l}(E_{(\bI,\bJ)})}
    {e_T(\mathcal{T}_{\bI}^-) e_T(\mathcal{T}_{\bJ}^+)}.
\end{align*}
\end{proof}
\end{lemma}


\section{Geometric Heisenberg and Clifford operators} \label{sec:geometric-operators}

\subsection{A complex of tautalogical bundles}

One has the \emph{tautological bundles}
\[
    V \times_{GL(V)} M(r,n) \to \M(r,n),\quad W \times \M(r,n) \to \M(r,n)
\]
where $M(r,n) = \{(A,B,i,j)\ |\ [A,B]+ij=0,\ (A,B,i,j) \text{ is
stable}\}$.  We denote these vector bundles by $V$ and $W$
respectively.  They are $T$-equivariant via the natural action of
$T$ on $M(r,n)$ and $W$.  We can then consider $A$ and $B$ to be
sections of the bundle $\Hom(V,V)$ and $i$ and $j$ to be sections of
the bundles $\Hom(W,V)$ and $\Hom(V,W)$ respectively. Over the
product $\M_\bc(r,n_1) \times \M_\mathbf{d}(r,n_2)$, we then have
the tautological bundles $V_k$ and $W_k$, and sections $A_k$, $B_k$,
$i_k$ and $j_k$, where $k=1,2$, coming from the tautological bundles
and sections on the $k$th factor.

We define a $T$-equivariant complex of vector bundles on
$\M_\bc(r,n_1)\times \M_{\mathbf{d}}(r,n_2)$ by
\begin{equation} \label{eq:master-complex}
   \Hom (V_1,V_2) \xrightarrow{\sigma}
        \begin{matrix} t\Hom(V_1, V_2) \oplus t^{-1}\Hom(V_1, V_2) \\
                        \oplus \\
                       \Hom(W_1, V_2) \oplus \Hom(V_1, W_2)
                       \end{matrix}
   {\xrightarrow{\tau}} \Hom (V_1,V_2),
\end{equation}
where $\sigma$ and $\tau$ are defined by
\begin{equation*}
   \sigma(\xi) = \begin{pmatrix} \xi A_1 - A_2 \xi \\
                            \xi B_1 - B_2 \xi \\
                            \xi i_1 \\
                            - j_2\xi \end{pmatrix}, \quad
   \tau \begin{pmatrix} C \\ D \\ I \\ J \end{pmatrix}
        = \left([A,D] + [C,B] + i_2J + Ij_1\right),
\end{equation*}
where
\[
    [A,D] = A_2 D - D A_1,\quad [C,B] = CB_1 - B_2C.
\]
One easily checks that $\tau \sigma=0$.

\begin{remark}
In the case when $\bc = \mathbf{d}$ and $n_1 = n_2$, the cohomology
of this complex is the tangent bundle to $\M_\bc(r,n_1)$ (more
precisely, the tangent bundle to the diagonal copy of
$\M_\bc(r,n_1)$ inside $\M_\bc(r,n_1) \times \M_\bc(r,n_2)$).  See
\cite{NY05}, in particular, the proof of Theorem~2.11.  Thus,
\eqref{eq:master-complex} can be seen be as a generalization of the
complex computing the tangent bundle.  It is also related to the
complex \cite[Equation~(5.1)]{Nak98} used by Nakajima to define
Hecke correspondences yielding the action of Kac-Moody algebras on
the homology of quiver varieties.  Nakajima identified a section of
the cohomology of that complex whose zero set defined the
correspondences used in his construction.
\end{remark}

\begin{lemma}  The cohomology $\ker \tau/\im \sigma$ is a vector bundle on
$\M_{\bc}(r,n_1) \times \M_{\mathbf{d}}(r,n_2)$.
\end{lemma}
\begin{proof}
We must show that $\tau$ is surjective and that $\sigma$ is
injective.  The argument is analogous to the proof of
\cite[Lemma~3.10]{Nak98} but since our notation and stability
condition differ from those used in \cite{Nak98}, we include the
proof here.

To show that $\sigma$ is injective, suppose $\xi \in \ker \sigma$.
Then $\ker \xi \subseteq V_1$ satisfies
\[
    \im i_1 \subseteq \ker \xi,\quad  A_1(\ker \xi) \subseteq \ker
    \xi,\quad  B_1(\ker \xi) \subseteq \ker \xi.
\]
Thus, stability of the point $(A_1,B_1,i_1,j_1)$ implies $\ker \xi =
V_1$ and so $\xi=0$.  Therefore $\sigma$ is injective.

To show that $\tau$ is surjective, suppose that $\zeta\in
\Hom(V_2,V_1)$ is orthogonal to $\im \tau$ with respect to the
(non-degenerate) trace pairing
\[
    \langle \cdot, \cdot \rangle : \Hom(V_1,V_2) \times \Hom(V_2,V_1) \longrightarrow
    \C.
\]
Then we have
\[
    A_1\zeta = \zeta A_2,\quad B_1\zeta = \zeta B_2,\quad \zeta i_2 = 0.
\]
Thus $\ker \zeta \subseteq V_2$ satisfies
\[
    \im i_2 \subseteq \ker \zeta,\quad  A_2(\ker \zeta) \subseteq \ker
    \zeta,\quad  B_2(\ker \zeta) \subseteq \ker \zeta,
\]
and so the stability condition for $(A_2,B_2,i_2,j_2)$ implies that
$\ker \zeta = V_2$.  Thus $\zeta = 0$ and therefore $\tau$ is
surjective.
\end{proof}

Denote the vector bundle $\ker \tau/\im \sigma$ by $\K_{\bc,
\mathbf{d}}(r,n_1,n_2)$.  The rank of $\K_{\bc,
\mathbf{d}}(r,n_1,n_2)$ is equal to $r(n_1 + n_2)$ as can be seen
from the following lemma.

\begin{lemma} \label{lem:master-complex-chern-class}
Let $(\bI,\bJ) \in \M_{\bc}(r,n_1)^T \times
\M_{\mathbf{d}}(r,n_2)^T$. The equivariant Euler class of the
restriction of $\K_{\bc, \mathbf{d}}(r,n_1,n_2)$ to the point
$(\bI,\bJ)$ is given as follows:
\begin{multline} \label{eq:master-complex-chern-class}
    e_T(\K_{\bc, \mathbf{d}}(r,n_1,n_2)_{(\bI,\bJ)}) =
    \prod_{\alpha,\beta = 1}^r \prod_{s \in \lambda(I^\alpha)} (b_\beta
    - b_\alpha + (d^\beta - c^\alpha - h_{I^\alpha,
    J^\beta}(s))\epsilon) \times \\
    \times \prod_{s\in \lambda(J^\alpha)}
    (b_\alpha - b_\beta +
    (d^\alpha - c^\beta + h_{J^\alpha, I^\beta}(s))\epsilon).
\end{multline}
\end{lemma}
\begin{proof}
The $T$-module decomposition of $\K_{\bc, \mathbf{d}}(r,n_1,n_2)$
for the case $\bc = \mathbf{d}= \mathbf{0}$ is computed in
\cite[Theorem~2.11]{NY05} (one must set $t_1 = t$ and $t_2 = t^{-1}$
there). The general result follows after replacing $b_\alpha$ by
$b_\alpha + c^\alpha \epsilon$ or $b_\alpha + d^\alpha \epsilon$ as
appropriate for all $\alpha \in \{1, \dots, r\}$.
\end{proof}


\subsection{Geometric Clifford operators}

Fix $l\in \{1,\hdots,r\}$.  We will see in
Corollary~\ref{cor:tnv-classes-clifford} that
\[
    c_{r(n_1 + n_2)}(\K_{\bc, \bc \pm 1_l}(r,n_1,n_2)) \ne 0,\quad n_1,
    n_2 \in \N,\ l \in \{1,\dots,r\},\ c \in \Z^r,
\]
and so we define
\[
    c_\tnv(\K_{\bc, \bc \pm 1_l}(r,n_1,n_2)) \stackrel{\text{def}}{=}
    c_{r(n_1+n_2)}(\K_{\bc, \bc \pm 1_l}(r,n_1,n_2)),
\]
the top non-vanishing equivariant Chern class of $\K_{\bc, \bc \pm
1_l}(r,n_1,n_2))$.

\begin{definition} \label{def:geom-cliff-operators}
For $l \in \{1,\hdots,r\}$ and $n\in \Z$, define operators
\[
    \Psi^l(n), \Psi^l(n)^* : \bigoplus_{\bc, k}
    \mathcal{H}^{2rk}_T(\mathcal{M}_\bc(r,k)) \to
    \bigoplus_{\bc, k} \mathcal{H}^{2rk}_T(\mathcal{M}_\bc(r,k))
\]
by
\begin{align*}
    \Psi^l(n)|_{\mathcal{H}^{2rk}_T(\M_\bc(r,k))} &=
    c_\tnv(\K_{\bc, \bc + 1_l}(r,k,k+n-c^l-1))
    \\
    &\in H^{2r(2k+n-c^l-1)}_T(\M_\bc(r,k)\times
    \M_{\bc+1_l}(r,k+n-c^l-1)), \\
    \Psi^l(n)^*|_{\mathcal{H}^{2rk}_T(\M_\bc(r,k))} &= c_\tnv(\K_{\bc,
    \bc - 1_l}(r,k,k-n+c^l)) \\
    &\in H^{2r(2k-n+c^l)} (\M_\bc(r,k) \times
    \M_{\bc-1_l}(r,k-n+c^l)).
\end{align*}
These operators will be called \emph{geometric Clifford operators}
(or \emph{geometric fermions}).
\end{definition}

\begin{lemma}
For $l \in \{1, \dots, r\}$ and $n \in \Z$, the operators
$\Psi^l(n)$ and $\Psi^l(n)^*$ are adjoint.
\end{lemma}
\begin{proof}
For $\bI \in \M_\bc(r,k)^T$, $\bJ \in \M_{\bc+1^l}(r,k+n-c^l-1)^T$,
we have
\begin{multline*}
    c_\tnv(\K_{\bc,\bc+1_l}(r,k,k+n-c^l-1)_{(\bI, \bJ)})\\
    = (-1)^{r(2k+n-c^l-1)}
    c_\tnv(\K_{\bc+1_l,\bc}(r,k+n-c^l-1,k)_{(\bJ, \bI)})
\end{multline*}
by Lemma~\ref{lem:master-complex-chern-class}.  Thus, by
Lemmas~\ref{lem:Euler-classes-tangent-spaces}
and~\ref{lem:structure-constants}, we have
\begin{align*}
    \langle \Psi^l(n) [\bI], [\bJ] \rangle &= \frac{
    c_\tnv(\K_{\bc,\bc+1_l}(r,k,k+n-c^l-1)_{(\bI, \bJ)})}
    {e_T(\mathcal{T}^-_\bI) e_T(\mathcal{T}_\bJ^+)} \\
    &= \frac{(-1)^{r(2k+n-c^l-1)} c_\tnv(\K_{\bc+1_l,\bc}(r,k+n-c^l-1,k)_{(\bJ, \bI)})}
    {(-1)^{rk} e_T(\mathcal{T}^+_\bI) (-1)^{r(k+n-c^l-1)} e_T(\mathcal{T}^-_\bJ)} \\
    &= \frac{
    c_\tnv(\K_{\bc+1_l,\bc}(r,k+n-c^l-1,k)_{(\bJ, \bI)})}{
    e_T(\mathcal{T}^+_\bI) e_T(\mathcal{T}^-_\bJ)} \\
    &= \langle \Psi^l(n)^* [\bJ], [\bI] \rangle.
\end{align*}
\end{proof}

\begin{theorem}\label{thm:clifford}
The geometric Clifford operators $\Psi^l(n)$, $\Psi^l(n)^*$ preserve
$\mathbf{A}$ and satisfy the relations
\begin{gather*}
    \{\Psi^l(n), \Psi^l(m)^*\} = \delta_{nm} \id,\quad
    \{\Psi^l(n), \Psi^l(m)\} = 0 = \{\Psi^l(n)^*, \Psi^l(m)^*\},\\
    [\Psi^l(n), \Psi^k(m)] = [\Psi^l(n), \Psi^k(m)^*] =
    [\Psi^l(n)^*, \Psi^k(m)^*] = 0,\quad l \ne k.
\end{gather*}
In particular, the maps
\[
    \psi^l(n) \mapsto \Psi^l(n),\quad \psi^l(n)^* \mapsto \Psi^l(n)^*,\quad
    n \in \Z,\,\quad l \in \{1, \dots, r\},
\]
define a representation of $\Cl$ on $\mathbf{A}$ and the linear map
$\mathbf{A} \to \mathbf{F}$ given by $[\bI] \mapsto \bI$ is an
isometric isomorphism of $\Cl$-modules. This isomorphism sends
$A_\bc(r,n)$ to $\mathbf{F}^{\bc}_n$.
\end{theorem}

The proof of this theorem, which involves a series of technical
combinatorial computations, will be given in
Section~\ref{sec:proof-clifford}.

\begin{corollary} \label{cor:tnv-classes-clifford}
For $n_1, n_2 \in \N$, $l \in \{1, \dots, r\}$,  and $\bc \in \Z^r$,
we have
\[
    c_{r(n_1+n_2)}(\K_{\bc, \bc \pm 1_l}(r,n_1,n_2)) \ne 0.
\]
\end{corollary}
\begin{proof}
Let $n=n_2-n_1+c^l+1$.  Then, by
Definition~\ref{def:geom-cliff-operators},
\[
    \Psi^l(n)|_{\mathcal{H}^{2rn_1}_T(\M_\bc(r,n_1))} =
    c_{r(n_1+n_2)}(\K_{\bc, \bc+1_l}(r,n_1,n_2)).
\]
We claim that $\psi^l(n)$ is a non-zero operator on
$\mathbf{F}^\bc_{n_1}$. If $n_1 = 0$, let $I^l = \vac{c^l}$.
Otherwise, if $n_2 > n_1$, let
\[
    I^l = (c^l+1) \wedge c^l \wedge (c^l-1) \wedge \dots \wedge
    (c^l-n_1+2) \wedge (c^l-n_1) \wedge (c^l-n_1-1) \wedge \dots,
\]
and if $n_2 \le n_1$, let
\begin{multline*}
    I^l = (c^l+1+n_2) \wedge c^l \wedge (c^l - 1) \wedge \dots\\
    \dots \wedge (c^l-n_1+n_2+2)
    \wedge (c^l-n_1 + n_2) \wedge (c^l-n_1+n_2-1) \wedge \dots.
\end{multline*}
One easily checks that $c(I^l) = c^l$ and $|I^l|=n_1$.  Thus if we
set
\[
    \bI = \left(\vac{c^1}, \dots, \vac{c^{l-1}}, I^l, \vac{c^{l+1}},
    \dots, \vac{c^r}\right),
\]
we have $\bc(\bI) = \bc$ and $|\bI|=n_1$. Since $\psi^l(n) (\bI) \ne
0$, we see that $\psi^l(n)$ is a non-zero operator on
$\mathbf{F}^\bc_{n_1}$. Therefore $\Psi^l(n)$ is a non-zero operator
on $\mathcal{H}^{2rn_1}_T(\M_\bc(r,n_1))$ by
Theorem~\ref{thm:clifford} and the first statement follows.  The
proof of the second statement is analogous.
\end{proof}
As mentioned above, this corollary justifies our calling
\[
  c_{r(n_1+n_2)} (\K_{\bc, \bc \pm 1_l}(r,n_1,n_2))
\]
the top non-vanishing equivariant Chern class of $\K_{\bc, \bc \pm
1_l}(r,n_1,n_2)$.


\subsection{Geometric Heisenberg operators} \label{sec:heisenberg-operators}

Let $(\bI, \bJ) \in \M_\bc(r,n_1)^T \times \M_\bc(r,n_2)^T$.  By
Lemma~\ref{lem:master-complex-chern-class}, we see that
\begin{multline*}
    e_T(\K_{\bc,\bc}(r,n_1,n_2)_{\bI,\bJ}) = \prod_{\alpha,
    \beta=1}^r \prod_{s \in \lambda(I^\alpha)} (b_\beta - b_\alpha +(c^\beta - c^\alpha -
    h_{I^\alpha, J^\beta}(s))\epsilon) \times \\
    \times \prod_{s \in \lambda(J^\alpha)} (b_\alpha - b_\beta +
    (c^\alpha - c^\beta  + h_{J^\alpha, I^\beta}(s)) \epsilon).
\end{multline*}
Now, if $n_1 \ne n_2$, then we must have $I^\alpha \ne J^\alpha$ for
some $\alpha$.  But then, since $c(I^\alpha)=c(J^\alpha)$, we have
$\lambda(I^\alpha) \ne \lambda(J^\alpha)$ and so $h_{I^\alpha,
J^\alpha}=0$ or $h_{J^\alpha, I^\alpha}=0$ by
Lemma~\ref{lem:hook=0}.  Therefore
\[
    e_T(\K_{\bc,\bc}(r,n_1,n_2)) = c_{r(n_1+n_2)}(\K_{\bc,\bc}(r,n_1,n_2))=0
    \quad \text{for} \quad n_1 \ne n_2
\]
by the Localization Theorem. We will see in
Corollary~\ref{cor:tnv-classes-heisenberg} that
\[
    c_{r(n_1+n_2)-1}(\K_{\bc,\bc}(r,n_1,n_2)) \ne 0
\]
and so we define
\[
    c_\tnv(\K_{\bc,\bc}(r,n_1,n_2)) \stackrel{\text{def}}{=}
    c_{r(n_1+n_2)-1}(\K_{\bc,\bc}(r,n_1,n_2)),
\]
the top non-vanishing equivariant Chern class of
$\K_{\bc,\bc}(r,n_1,n_2)$.

\begin{lemma}\label{lem:bosonic-chern-class}
Let $(\bI,\bJ) \in \M_\bc(r,n_1)^T \times \M_\bc(r,n_2)^T$ such that
$\lambda(I^l) \ne \lambda(J^l)$ and let $k$ be the smallest integer
such that $\lambda(I^l)_k \ne \lambda(J^l)_k$. We have
\begin{multline*}
    c_\tnv(\K_{\bc,\bc}(r,n_1,n_2)_{(\bI, \bJ)})=
    \prod_{\alpha,\beta=1}^r  \prod_{\begin{matrix} {\scriptstyle
    s \in \lambda(I^\alpha),}\\ {\scriptstyle s \ne
    (k,\lambda(I^l)_k)} \\ {\scriptstyle
    \text{ if } \alpha=\beta=l} \end{matrix}} (b_ \beta - b_\alpha
    +(c^\beta - c^\alpha - h_{I^\alpha, J^\beta}(s)) \epsilon) \times \\
    \times \prod_{s \in \lambda(J^\alpha)}
    (b_\alpha - b_\beta + (c^\alpha - c^\beta + h_{J^\alpha, I^\beta}(s)) \epsilon)
\end{multline*}
if $\lambda(I^l)_k > \lambda(J^l)_k$, and
\begin{multline*}
    c_\tnv(\K_{\bc,\bc}(r,n_1,n_2)_{(\bI,\bJ)})=
    \prod_{\alpha,\beta=1}^r \prod_{s \in \lambda(I^\alpha)} (b_\beta -
    b_\alpha +(c^\beta - c^\alpha -h_{I^\alpha, J^\beta}(s))
    \epsilon) \times \\
    \times \prod_{\begin{matrix} {\scriptstyle
    s \in \lambda(J^\alpha),}\\ {\scriptstyle s \ne (k,\lambda(J^l)_k)
    \text{ if } \alpha=\beta=l} \end{matrix}} (b_\alpha -
    b_\beta + (c^\alpha - c^\beta + h_{J^\alpha, I^\beta}(s)) \epsilon).
\end{multline*}
if $\lambda(J^l)_k > \lambda(I^l)_k$.
\end{lemma}

\begin{proof}
By the above comments, $c_\tnv(\K_{\bc,\bc}(r,n_1,n_2)_{(\bI,\bJ)})$
is obtained by removing one factor of zero from the product
appearing in \eqref{eq:master-complex-chern-class}.  The result then
follows from Lemma~\ref{lem:hook=0}.
\end{proof}

\begin{remark}
If $(\bI,\bJ) \in \M_\bc(r,n_1)^T \times \M_\bc(r,n_2)^T$ such that
$\lambda(I^l) \ne \lambda(J^l)$ for more than one choice of $l$,
then $c_\tnv(\K_{\bc,\bc}(r,n_1,n_2)_{\bI, \bJ})=0$. Thus it does
not matter which $l$ we choose in
Lemma~\ref{lem:bosonic-chern-class}.
\end{remark}

If $n_1 = n_2$, then we will see in
Corollary~\ref{cor:tnv-classes-heisenberg} that the $r(n_1+n_2)$-th
equivariant Chern class of $\K_{\bc,\bc}(r,n_1,n_2)$ does \emph{not}
vanish and so we define
\[
    c_\tnv(\K_{\bc,\bc}(r,n,n)) \stackrel{\text{def}}{=}
    c_{2rn}(\K_{\bc,\bc}(r,n,n)).
\]

\begin{lemma} \label{lem:bosonic-chern-class-diagonal}
Let $(\bI, \bJ) \in \M_\bc(r,n)^T \times \M_\bc(r,n)^T$.  Then
\[
    c_\tnv(\K_{\bc,\bc}(r,n,n)_{\bI, \bJ}) = 0
\]
if $\bI \ne \bJ$ and
\[
    c_\tnv(\K_{\bc,\bc}(r,n,n)_{\bI,\bI}) = e_T(\mathcal{T}_\bI).
\]
\end{lemma}
\begin{proof}
The case $\bI \ne \bJ$ follows as above.  The case $\bI = \bJ$
follows immediately from
Lemmas~\ref{lem:Euler-classes-tangent-spaces}
and~\ref{lem:master-complex-chern-class}.
\end{proof}

Fix $n \neq m$ and consider the action of the subtorus $(\C^*)^r
\subset T$ on $\M_\bc(r,n)\times \M_\bc(r,m)$, acting on the framing
over $l_\infty \subset \mathbb{P}^2$.  The connected components of
the fixed point set of this action are products of Hilbert schemes
\begin{multline*}
    \left(\M_\bc(r,n)\times \M_\bc(r,m)\right)^{(\C^*)^r} \\
    = \bigsqcup_{\mathbf{n}, \mathbf{m}}
    \left({\C^2}^{[n^1]}\times \dots \times {\C^2}^{[n^r]}\right)
    \times \left({\C^2}^{[m^1]}\times \dots \times {\C^2}^{[m^r]}\right),
\end{multline*}
where the union is over $\mathbf{n} = (n^1,\dots,n^r)$, $\mathbf{m}=
(m^1,\dots,m^r)$ with $\sum_i n^i = n$, $\sum_i m^i = m$, and
${\C^2}^{[k]}$ denotes the Hilbert scheme of $k$ points in $\C^2$.
We have the inclusions
\begin{multline*}
    i_{\mathbf{n},\mathbf{m}}: \left({\C^2}^{[n^1]}\times \dots \times {\C^2}^{[n^r]}\right)
    \times \left({\C^2}^{[m^1]}\times \dots \times
    {\C^2}^{[m^r]}\right) \\
    \hookrightarrow \M_\bc(r,n)\times \M_\bc(r,m).
\end{multline*}
Define disjoint sets $A_l$, $l \in \{1,\dots,r\}$, by
\[
    A_l = \{(\mathbf{n},\mathbf{m})\ |\ n^\alpha = m^\alpha \ \text{for}\ \alpha<l,
    \ n^l \neq m^l \}.
\]
Then, for $l \in \{1,\dots,r\}$, $n \ne m$, we have a partition of
the $(\C^*)^r$ fixed-point components
\begin{gather*}
    \left(\M_\bc(r,n)\times \M_\bc(r,m)\right)^{(\C^*)^r} = \bigsqcup_{l=1}^r
    \mathcal{A}_l,\\
    \mathcal{A}_l = \bigsqcup_{(\mathbf{n},\mathbf{m}) \in A_l}
     \left({\C^2}^{[n^1]}\times \dots \times {\C^2}^{[n^r]}\right)
    \times \left({\C^2}^{[m^1]}\times \dots \times {\C^2}^{[m^r]}\right),
\end{gather*}
and the associated inclusions $i_l : \mathcal{A}_l \hookrightarrow
\M_\bc(r,n) \times \M_\bc(r,m)$.  By the Localization Theorem, the
restriction
\begin{multline*}
    i^* = \sum_{\mathbf{n},\mathbf{m}} i_{\mathbf{n},\mathbf{m}}^* :
    \mathcal{H}^*_T \left(\M_\bc(r,n)\times \M_\bc(r,m)\right) \\
    \longrightarrow
    \mathcal{H}^*_T \left((\M_\bc(r,n)\times \M_\bc(r,m))^{(\C^*)^r}\right)
\end{multline*}
is an isomorphism.  By construction, the classes
\[
    \tilde \gamma^l \stackrel{\text{def}}{=} i_l^*(1) = \sum_{(\mathbf{n}, \mathbf{m})
    \in A_l} i_{\mathbf{n}, \mathbf{m}}^*(1) \in \mathcal{H}^0_T(\mathcal{A}_l)
    \subseteq
    \mathcal{H}^0_T\left((\M_\bc(r,n)\times \M_\bc(r,m))^{(\C^*)^r}\right)
\]
are orthogonal idempotents ($\tilde \gamma^k \cup \tilde \gamma^l =
\delta_{k,l} \tilde \gamma^k$), and they decompose the identity:
\[
    1 = \sum_{l=1}^r \tilde \gamma^l \in  \mathcal{H}^0_T\left((\M_\bc(r,n)\times
    \M_\bc(r,m))^{(\C^*)^r}\right).
\]
For $l \in \{1,\dots,r\}$, let
\[
    \gamma^l = \epsilon \cup (i^*)^{-1}(\tilde \gamma^l) \in \mathcal{H}^2_T( (\M_\bc(r,n)
    \times \M_\bc(r,m)).
\]
It follows from the definitions that
\[
    \gamma^l_{\bI, \bJ} = \begin{cases} \epsilon & \text{if } |I^\alpha| =
    |J^\alpha| \text{ for } \alpha < l \text{ and } |I^l|\ne|J^l|, \\
    0 & \text{otherwise}, \end{cases}
\]
where $\gamma^l_{\bI, \bJ} = i^*_{\bI, \bJ}(\gamma^l)$ and $i_{\bI,
\bJ} : (\bI, \bJ) \hookrightarrow \M_\bc(r,n) \times \M_\bc(r,m)$ is
the inclusion of the fixed point.

\begin{definition} \label{def:geom-heis-operators}
For $l \in \{1,\hdots,r\}$ and $n\in \Z$, define an operator
\[
    P^l(n) : \bigoplus_{\bc, k} \mathcal{H}^{2rk}_T(\mathcal{M}_\bc(r,k)) \to
    \bigoplus_{\bc, k} \mathcal{H}^{2rk}_T(\mathcal{M}_\bc(r,k))
\]
by
\begin{align*}
    P^l(n)|_{\mathcal{H}^{2rk}_T(\M_\bc(r,k))} &= \begin{cases}
    \gamma^l \cup c_\tnv(\K_{\bc,\bc}(r,k,k-n)) & n < 0,\\
    -\gamma^l \cup c_\tnv (\K_{\bc,\bc}(r,k,k-n)) & n > 0,
    \end{cases} \\
    &\qquad \in \H^{2r(2k-n)}_T(\M_\bc(r,k)\times \M_\bc(r,k-n)), \\
    P^l(0)|_{\mathcal{H}^{2rk}_T(\M_\bc(r,k))} &= c^l c_\tnv(\K_{\bc,\bc}(r,k,k))
    = c^l \id
\end{align*}
(the last equality follows from Lemmas~\ref{lem:structure-constants}
and~\ref{lem:bosonic-chern-class-diagonal}). These operators will be
called \emph{geometric Heisenberg operators} (or \emph{geometric
bosons)}.
\end{definition}

\begin{remark}
We motivate the presence of the classes $\gamma^l$ in
Definition~\ref{def:geom-heis-operators}.  Note that the
decomposition $\epsilon = \sum_{l=1}^r \gamma^l$ implies
\[
    \epsilon \cup c_\tnv(\K_{\bc,\bc}(r,k,k-n)) =
    \sum_{l=1}^r P^l(n)|_{\mathcal{H}^{2rk}_T(\M_\bc(r,k))}.
\]
Thus, we have decomposed the operator $\epsilon \cup
c_\tnv(\K_{\bc,\bc}(r,k,k-n))$ as a sum of $r$ different operators.
In particular, for $r=1$, the Heisenberg operators are simply given
by the classes $\epsilon \cup c_\tnv(\K_{\bc,\bc}(r,k,k-n))$.
\end{remark}

\begin{lemma}
For $n \in \Z$, $l \in \{1, \dots, n\}$, the operators $P^l(n)$ and
$P^l(-n)$ are adjoint.
\end{lemma}
\begin{proof}
For $n = 0$, the statement is obvious.  Assume $n>0$.  For $\bI \in
\M_\bc(r,k)^T$, $\bJ \in \M_\bc(r,k-n)^T$, we have
\[
    c_\tnv(\K_{\bc,\bc}(r,k,k-n)_{(\bI, \bJ)}) = (-1)^{r(2k-n)-1}
    c_\tnv(\K_{\bc,\bc}(r,k-n,k)_{(\bJ, \bI)})
\]
by Lemma~\ref{lem:bosonic-chern-class}.  Thus, by
Lemmas~\ref{lem:Euler-classes-tangent-spaces}
and~\ref{lem:structure-constants}, we have
\begin{align*}
    \langle P^l(n) [\bI], [\bJ] \rangle &= \frac{-\gamma^l_{\bI,\bJ} \cup
    c_\tnv(\K_{\bc,\bc}(r,k,k-n)_{(\bI,
    \bJ)})}{e_T(\mathcal{T}^-_\bI) e_T(\mathcal{T}_\bJ^+)} \\
    &= \frac{-(-1)^{r(2k-n)-1} \gamma^l_{\bI,\bJ} \cup
    c_\tnv(\K_{\bc,\bc}(r,k-n,k)_{(\bJ, \bI)})}{(-1)^{rk}
    e_T(\mathcal{T}^+_\bI) (-1)^{r(k-n)} e_T(\mathcal{T}^-_\bJ)} \\
    &= \frac{\gamma^l_{\bI,\bJ} \cup
    c_\tnv(\K_{\bc,\bc}(r,k-n,k)_{(\bJ, \bI)})}{
    e_T(\mathcal{T}^+_\bI) e_T(\mathcal{T}^-_\bJ)} \\
    &= \langle P^l(-n) [\bJ], [\bI] \rangle.
\end{align*}
\end{proof}

\begin{theorem} \label{thm:Heisenberg}
The geometric Heisenberg operators $P^l(n)$ preserve $\mathbf{A}$
and satisfy the relations
\[
    [P^k(n), P^l(0)]=0,\quad [P^k(m),P^l(n)] = \frac{1}{m}\delta_{m,-n}\delta_{k,l}
    \id,\ m \ne 0.
\]
In particular, the maps
\[
    p^l(n) \mapsto P^l(n),\ n \in \Z,\, l \in \{1, \dots, r\},\quad K \mapsto
    \id,
\]
define a representation of $\mathfrak{s}$ on $\mathbf{A}$ and the
linear map $\mathbf{A} \to \mathbf{B}$ given by
\[
[\bI] \mapsto (q^{c(I^1)} s_{\lambda(I^1)}, \dots, q^{c(I^r)}
s_{\lambda(I^r)})
\]
is an isometric isomorphism of $\mathfrak{s}$-modules.  This
isomorphism identifies $A_\bc(r,n)$ with $\mathbf{B}^{\bc}_n$.
\end{theorem}
The proof of this theorem, which involves a series of technical
combinatorial computations, will be given in
Section~\ref{sec:proof-heisenberg}.

\begin{corollary} \label{cor:tnv-classes-heisenberg}
For $n_1, n_2 \in \N$, $n_1 \ne n_2$, and $\bc \in \Z^r$, we have
\[
    c_{r(n_1+n_2)-1} (\K_{\bc, \bc}(r,n_1,n_2)) \ne 0 \quad
    \text{and} \quad
    c_{2rn_1}(\K_{\bc,\bc}(r,n_1,n_1)) \ne 0.
\]
\end{corollary}
\begin{proof}
Let $n=n_1-n_2$.  Then, by Definition~\ref{def:geom-heis-operators},
\[
    P^l(n)|_{\mathcal{H}^{2rn_1}_T(\M_\bc(r,n_1))} = \pm \gamma^l \cup
    c_{r(n_1+n_2)-1} (\K_{\bc, \bc}(r,n_1,n_2)).
\]
It is easily seen from the description of the operators $p^l(n)$ in
Section~\ref{sec:bosonic-fock-space} that $p^l(n)$ is a non-zero
operator on $B^\bc_{n_1}$ (note that our assumptions imply $n \le
n_1$).  Therefore, $P^l(n)$ is a non-zero operator on
$\mathcal{H}^{2rn_1}_T(\M_\bc(r,n_1))$ by
Theorem~\ref{thm:Heisenberg} and the first result follows.  The
second follows analogously from the fact that $p^l(0)$ is a non-zero
operator.
\end{proof}
As mentioned above, this corollary justifies our calling
\[
  c_{r(n_1+n_2)-1} (\K_{\bc, \bc}(r,n_1,n_2)) \quad \text{and} \quad
  c_{2rn} (\K_{\bc, \bc}(r,n,n))
\]
the top non-vanishing equivariant Chern classes of $\K_{\bc,
\bc}(r,n_1,n_2)$, $n_1 \ne n_2$, and $\K_{\bc, \bc}(r,n,n)$
respectively.


\subsection{Sheaf theoretic interpretation of the geometric Heisenberg
and Clifford operators} \label{sec:sheaf-interpretation}

Using arguments similar to those appearing in
\cite[Section~2]{NY05}, one can interpret our geometric Heisenberg
and Clifford operators in the language of sheaves.  In so doing, we
see that in the rank 1 case they are closely related to operators
defined by Carsslon-Okouknov in \cite{CO08}.  In particular, we have
the following.

\begin{proposition}
Let $((E_1, \Phi_1), (E_2, \Phi_2)) \in \M_\bc(r,n_1) \times
\M_\mathbf{d}(r,n_2)$ be a pair of framed torsion-free sheaves. Then
\[
  \K_{\bc,\mathbf{d}}(r,n_1,n_2)_{((E_1, \Phi_1), (E_2, \Phi_2))}
  = \Ext^1(E_1,E_2(-l_\infty)).
\]

\begin{proof}
Let $(E, \Phi) \in \M_\bc(r,n)$ be a framed torsion-free sheaf with
corresponding ADHM data $(A,B,i,j)$ (see Section~\ref{sec:Mrn}).
The tangent space to $\M_\bc(r,n)$ at $(E, \Phi)$ can be computed in
two ways, one using sheaf theory and the other using the ADHM
description (see \cite{NY05}, in particular, the proof of
Theorem~2.11).  In the first case, the tangent space at $(E, \Phi)$
is given by $\Ext^1(E,E(-l_\infty))$. In the second, the tangent
space at $(E, \Phi)$ is the middle cohomology of the complex
\begin{equation} \label{eq:tangent-complex}
   \Hom (V,V) \xrightarrow{\sigma}
        \begin{matrix} t\Hom(V, V) \oplus t^{-1}\Hom(V, V) \\
                        \oplus \\
                       \Hom(W, V) \oplus \Hom(V, W)
                       \end{matrix}
   {\xrightarrow{\tau}} \Hom (V,V),
\end{equation}
where $\sigma$ and $\tau$ are defined by
\begin{equation*}
   \sigma(\xi) = \begin{pmatrix} \xi A - A \xi \\
                            \xi B - B \xi \\
                            \xi i \\
                            - j\xi \end{pmatrix}, \quad
   \tau \begin{pmatrix} C \\ D \\ I \\ J \end{pmatrix}
        = [A,D] + [C,B] + iJ + Ij.
\end{equation*}

Now consider the particular case of the framed rank $2r$
torsion-free sheaf
\[
    (E, \Phi) = (E_1\oplus E_2, \Phi_1 \oplus \Phi_2)
\]
which is a point of $\M_{(\bc,\mathbf{d})}(2r,n_1+n_2)$ where $(\bc,
\mathbf{d}) = (c^1, \dots, c^r, d^1, \dots, d^r)$.  In the language
of sheaves, the tangent space at $(E, \Phi)$ is then
\[
    \Ext^1(E_1\oplus E_2, E_1(-l_\infty)\oplus E_2(-l_\infty))
\]
while in the ADHM description, it is the cohomology of the complex
\begin{multline} \label{eq:tangent-complex2}
   \Hom (V_1\oplus V_2,V_1\oplus V_2) \\ \xrightarrow{\sigma}
        \begin{matrix} t\Hom(V_1\oplus V_2, V_1\oplus V_2) \oplus
        t^{-1}\Hom(V_1\oplus V_2, V_1\oplus V_2) \\
                        \oplus \\
                       \Hom(W_1\oplus W_2, V_1\oplus V_2) \oplus
                       \Hom(V_1\oplus V_2, W_1\oplus W_2)
                       \end{matrix} \\
   {\xrightarrow{\tau}} \Hom (V_1\oplus V_2,V_1\oplus V_2).
\end{multline}

Define a $\C^*$-action on $\M_{(\bc,\mathbf{d})}(2r,n_1+n_2)$ using
the one-parameter subgroup
\[
    s \mapsto \id_{W_1} \oplus s \id_{W_2} \in GL(W_1)\times GL(W_2)
    \subset GL(W_1 \oplus W_2).
\]
The fixed points of this $\C^*$-action are those rank $2r$
torsion-free sheaves which are isomorphic to a direct sum of two
rank $r$ torsion-free sheaves (see \cite[Proposition~2.9]{NY05}). In
particular, $(E_1\oplus E_2, \Phi_1 \oplus \Phi_2)$ is fixed by this
action and thus the tangent space $\Ext^1(E_1\oplus E_2,
E_1(-l_\infty)\oplus E_2(-l_\infty))$ has an induced $\C^*$-action
and decomposes into isotypic components for this action. Over the
$\C^*$-fixed point $(E_1 \oplus E_2, \Phi_1 \oplus \Phi_2)$, the
complex \eqref{eq:tangent-complex2} also decomposes into isotypic
components. As in the proof of \cite[Theorem~2.11]{NY05}, ones see
that the complex \eqref{eq:master-complex} is the isotypic
subcomplex of \eqref{eq:tangent-complex2} of weight 1. Similarly the
subbundle $\Ext^1(E_1,E_2(-l_\infty))$ is the isotypic subbundle of
$\Ext^1 (E_1 \oplus E_2, E_1(-l_\infty) \oplus E_2(-l_\infty))$ of
weight 1.
\end{proof}
\end{proposition}

Consider now the case $r=1$ and recall that $\M(1,n)$ is the Hilbert
scheme of $n$ points in $\C^2$.  Note that if one takes $\mathcal{L}
= \mathcal{O}(-l_\infty)$ in \cite[Section~1.2]{CO08}, then
$\chi(\mathcal{L})=0$ and thus
\[
    \mathrm{E}_{(I,J)} = -\chi(I, J(-l_\infty)) = - \sum_{i=0}^2
    (-1)^i \Ext^i(I, J(-l_\infty)) = \Ext^1(I, J(-l_\infty))
\]
for ideal sheaves $I$ and $J$.  Here we use the fact that $\Ext^0(I,
J(-l_\infty)) = \Ext^2(I, J(-l_\infty))=0$ (see
\cite[Proposition~2.1]{NY05}).  Therefore, in the rank one case, our
vector bundle $\K_{\bc, \mathbf{d}}(1,n_1,n_2)$ is an example of the
virtual vector bundles considered in \cite{CO08} with a modified
torus action.


\section{Proof of Theorem~\ref{thm:clifford}}  \label{sec:proof-clifford}
\subsection{Notation}

In order to simplify notation in what follows, for $r$-colored
semi-infinite monomials $\bI = (I^1,\hdots,I^r)$ and $\bJ =
(J^1,\hdots,J^r)$ with $\bc(\bI) + 1_l = \bc(\bJ)$ for some $l \in
\{1,\dots,r\}$, define
\[
    f_{\bI,\bJ} =
    c_\tnv \left(\K_{\bc(\bI), \bc(\bJ)}(r, |\bI|, |\bJ|)_{(\bI,\bJ)} \right).
\]
When non-zero, $f_{\bI,\bJ}$ is homogeneous of degree $2r(|\bI| +
|\bJ|)$. When $r=1$ and $c(J) = c(I)+1$,
\[
    f_{I,J} = (-1)^{|I|} \epsilon^{|I| + |J|} \prod_{s \in \lambda(I)}
    (a_I(s) + l_J(s)) \prod_{s\in \lambda(J)}
    (a_J(s) + l_I(s) + 2).
\]
For $\bI \in \M_\bc(r,n_1)^T$, $\bJ \in \M_{\mathbf{d}}(r,n_2)^T$,
we also define the polynomial
\[
    d_{\bI,\bJ} = e_T(\mathcal{T}_\bI^-)e_T(\mathcal{T}_\bJ^+)
    \in \C[b_1,\hdots,b_r,\epsilon].
\]
Note that when $r=1$, we have
\[
    d_{I, J} = (-1)^{|I|} h_{I} h_{J} \epsilon^{|I| + |J|}.
\]
Thus
\[
    f_{I,J} = (-1)^{|I|}\tilde{f}_{I,J} \epsilon^{|I|
    + |J|}, \quad
    d_{I,J} = (-1)^{|I|}\tilde{d}_{I,J} \epsilon^{|I|
    + |J|},
\]
where
\[
    \tilde{f}_{I,J}=\prod_{s \in \lambda(I)} (a_I(s) + l_J(s))
\prod_{s\in \lambda(J)} ( a_J(s) + l_I(s) + 2), \quad
    \tilde{d}_{I,J}= h_I h_J.
\]


\subsection{Combinatorics of the rank one fermionic Fock space}
\label{sec:fermionic-combinatorics-rank-one}

Let $I = i_1 \wedge i_2 \wedge \hdots$ be a semi-infinite monomial
of charge $c$, and let $J = j_1 \wedge j_2 \wedge \hdots$ be a
semi-infinite monomial of charge $c + 1$.

\begin{lemma} \label{lem:i_1=j_1}
Suppose $i_1 = j_1$.  Let $l$ be the largest positive integer such
that $i_l = i_1 - l + 1$ (in other words, $\lambda(I)_l =
\lambda(I)_1$) and let
\[
    I' = i_{l+1} \wedge i_{l+2} \wedge \dots,\quad J' = j_{l+1}
    \wedge j_{l+2} \wedge \dots.
\]
Then we have
\[
    \tilde{f}_{I,J} = \begin{cases}
        (-1)^l \tilde{f}_{I',J'}\prod_{s \in R(I)} h_I(s) \prod_{s\in R(J)}
        h_J(s) & \text{if } i_k = j_k \
        \forall \ 1 \le k \le l,\\
        0 & \text{otherwise}, \end{cases}
\]
where
\begin{align*}
    R(I) &=
    (\lambda(I)_1, \lambda(I)_2, \dots, \lambda(I)_l),\quad \text{and} \\
    R(J) &=
    (\lambda(J)_1, \lambda(J)_2, \dots, \lambda(J)_l).
\end{align*}

\begin{proof}
Suppose $i_k \ne j_k$ for some $1 \le k \le l$.  Choose $k$ to be
minimal. Then $\lambda(J)_p = \lambda(I)_p - 1 = \lambda(I)_1-1$ for
all $1 \le p \le k-1$ and $\lambda(J)_k < \lambda(I)_k-1$.  Let
$s=(k,\lambda(I)_k-1)$.  Then
\[
    a_I(s) + l_J(s) = 1 + (-1) = 0
\]
and the result follows.  Therefore we assume $i_k = j_k$ for all $1
\le k \le l$. So for all $1 \le k \le l$ we have $\lambda(I)_k =
\lambda(I)_1$, $\lambda(J)_k = \lambda(J)_1$ and $\lambda(I)_k =
\lambda(J)_k +1$.  Also, by our choice of $l$, we have
$\lambda(I)_{l+1} < \lambda(I)_l = \lambda(I)_1$.  Note that
\begin{align*}
    \lambda(I') &= (\lambda(I)_{l+1}, \lambda(I)_{l+2}, \dots),\quad \text{and}\\
    \lambda(J') &= (\lambda(J)_{l+1}, \lambda(J)_{l+2}, \dots).
\end{align*}
We have
\begin{align*}
    \tilde{f}_{I,J} &= \prod_{s\in R(I)} (a_I(s) + l_J(s)) \prod_{s\in R(J)}
    (a_J(s) + l_I(s) +2) \\
    & \qquad \qquad \prod_{s\in\lambda(I')} (a_{I'}(s) + l_{J'}(s))
    \prod_{s\in\lambda(J')} (a_{J'}(s) + l_{I'}(s) +2) \\
    &= \tilde{f}_{I',J'}\prod_{s\in R(I)} (a_I(s) + l_J(s))
    \prod_{s\in R(J)}  (a_J(s) + l_I(s) +2).
\end{align*}
Now, for $s \in R(J)$, we have $a_I(s) = a_J(s) + 1$.  Therefore
\[
    \prod_{s\in R(I)} (a_I(s) + l_J(s)) = \prod_{s\in R(I) \setminus R(J)}
    (a_I(s) + l_J(s)) \prod_{s\in R(J)} h_J(s),
\]
and
\[
    \prod_{s\in R(J)}  (a_J(s) + l_I(s) +2) =
    \prod_{s\in R(J)} h_I(s).
\]
Now,
\[
    \prod_{s\in R(I) \setminus R(J)} (a_I(s) + l_J(s)) = (-1)(-2) \cdots (-l)
    = (-1)^l \prod_{s\in R(I) \setminus R(J)} h_I(s).
\]
Therefore
\[
    \tilde{f}_{I,J} = (-1)^l \tilde{f}_{I',J'}\prod_{s \in R(I)} h_I(s) \prod_{s\in R(J)}
    h_J(s).
\]
\end{proof}
\end{lemma}

\begin{lemma} \label{lem:i_k=j_k+1} If $i_k = j_{k+1}$ for all $k
\ge 1$, then
\[
    \tilde{f}_{I,J} = h_I h_J.
\]

\begin{proof}
The hypotheses imply that the partition $\lambda(I)$ is obtained
from the partition $\lambda(J)$ by removing the largest part
$\lambda(J)_1$. Therefore,
\[
    l_J(s) = l_I(s) + 1,\quad s \in \lambda(J),
\]
and so
\begin{align*}
    \prod_{s \in \lambda(I)} (a_I(s) + l_J(s)) &= \prod_{s \in \lambda(I)}
    (a_I(s) + l_I(s) + 1) = h_I,\\
    \prod_{s \in \lambda(J)} (a_J(s) + l_I(s) + 2)
    &= \prod_{s \in \lambda(J)} (a_J(s) + l_J(s) + 1) = h_J,
\end{align*}
and the result follows.
\end{proof}
\end{lemma}

\begin{lemma} \label{lem:0 Case}
If $i_1\neq j_1$ and $i_k \neq j_{k+1}$ for some $k \ge 1$, then
$\tilde{f}_{I,J} = 0$.

\begin{proof}
Suppose $i_1 \ne j_1$ and let $k$ be the smallest positive integer
such that $i_k \ne j_{k+1}$.  Thus $k$ is also the smallest positive
integer such that $\lambda(I)_k \ne \lambda(J)_{k+1}$. Suppose
$j_{k+1} > i_k$. Then $\lambda(J)_{k+1} > \lambda(I)_k$. Set
$s=(k+1,\lambda(J)_{k+1})$.  Then $a_J(s)=0$.  Also, since
$\lambda(I)_{k-1} = \lambda(J)_k \ge \lambda(J)_{k+1}$, we have
$l_I(s) = -2$. Then
\[
a_J(s) + l_I(s) + 2 = 0 + (-2) + 2 = 0.
\]
Now suppose $j_{k+1} < i_k$.  Then  $\lambda(J)_{k+1} <
\lambda(I)_k$.  If $s=(k,\lambda(I)_k)$, then $a_I(s)=0$.  Also, for
$k > 1$, since $\lambda(J)_k = \lambda(I)_{k-1} \ge \lambda(I)_k$,
we have $l_J(s) = 0$.  Therefore
\[
a_I(s) + l_J(s) = 0.
\]
If $k=1$, then either $\lambda(J)_1 \ge \lambda(I)_1$, in which case
the above still holds, or $\lambda(J)_1 \le \lambda(I)_1 - 2$ (we
cannot have $\lambda(J)_1 = \lambda(I)_1 -1$ since this would imply
$i_1 = j_1$).  In this case, let $s=(1,\lambda(I)_1-1)$.  Then
$a_I(s)=1$ and $l_J(s)=-1$ and we have
\[
a_I(s) + l_J(s) = 0.
\]
\end{proof}
\end{lemma}

\begin{proposition} \label{prop:ferm-poly-identity-rank-one}
We have
\[
    f_{I,J} = \begin{cases} (-1)^n d_{I,J} & \text{if } J = (-1)^n \psi(k) I
    \text{ for some } k, n \in \Z,\\
    0 & \text{otherwise}. \end{cases}
\]

\begin{proof}
It suffices to consider only the coefficients $\tilde f_{I,J}$ and
$\tilde d_{I,J}$.  We first note that there exists a $k$ such that
$J = \pm \psi(k)I$ if and only if there exists an $l$ such that $i_m
= j_m$ for all $m < l$, $j_l=k$ and $j_{m+1}=i_m$ for all $m \ge l$.
In this case, we have $J = (-1)^{l-1} \psi(k) I$.

Let $l$ be the smallest positive integer such that $i_l \ne j_l$ and
let
\[
    I' = i_l \wedge i_{l+1} \wedge \dots,\quad J' = j_l \wedge j_{l+1} \wedge \dots.
\]
Repeated application of Lemma~\ref{lem:i_1=j_1} gives
\[
    \tilde{f}_{I,J} = (-1)^{l-1} \tilde{f}_{I',J'} \prod_{s \in R(I)}
    h_I(s) \prod_{s \in R(J)} h_J(s).
\]
By Lemmas~\ref{lem:i_k=j_k+1} and~\ref{lem:0 Case},
$\tilde{f}_{I',J'}=0$ unless $i_m = j_{m+1}$ for all $m \ge l$, in
which case
\[
    \tilde{f}_{I',J'} = h_{I'} h_{J'},
\]
and so
\[
    \tilde{f}_{I,J} = (-1)^{l-1} h_I h_J = (-1)^{l-1} \tilde
    d_{I,J}.
\]
Thus, we have
\[
    \tilde{f}_{I,J} = \begin{cases} (-1)^n \tilde d_{I,J} &
    \text{if } J = (-1)^n \psi(k) I \text{ for some } k, n \in \Z,\\
    0 & \text{otherwise}. \end{cases}
\]
as desired.
\end{proof}
\end{proposition}


\subsection{Combinatorics of the $r$-colored fermionic Fock space}
\label{sec:fermionic-combinatorics-rank-r}

Suppose $\bI$ and $\bJ$ are $r$-colored semi-infinite monomials and
$\alpha, \beta \in \{1, \dots, r\}$. In order to simplify notation
in the following proofs, we define
\[
    X^{\alpha, \beta}_{\bI,\bJ} = \prod_{s \in \lambda(I^\alpha)}
    (b_\alpha - b_\beta + (c(I^\alpha)
    - c(J^\beta) + h_{I^\alpha,J^\beta}(s)) \epsilon)
\]

\begin{proposition} \label{prop:ferm-poly-identity-rank-r}
If $\bI$ and $\bJ$ are semi-infinite monomials with $\bc(\bI) + 1_l
= \bc(\bJ)$, then
\[
    f_{\bI, \bJ} = \begin{cases}
    (-1)^{n} d_{\bI, \bJ} & \text{if }
    \bJ = (-1)^n \psi^l(k) \bI \text{ for some } k,n \in \Z, \\
    0 & \text{otherwise}. \end{cases}
\]
\end{proposition}

\begin{proof}
We have
\begin{equation} \label{eq:fermionic-polynomial}
    f_{\bI, \bJ} =
    \prod_{\alpha,\beta = 1}^r (-1)^{|I^\alpha|} X^{\alpha, \beta}_{\bI, \bJ}
    X^{\alpha, \beta}_{\bJ, \bI} = f_{I^l, J^l}
    \prod_{\begin{array}{cc}
    {\scriptstyle 1 \le \alpha, \beta \le r} \\ {\scriptstyle
    \alpha \ne l \text{ or } \beta \ne l}
    \end{array}} (-1)^{|I^\alpha|} X^{\alpha, \beta}_{\bI, \bJ}
    X^{\alpha, \beta}_{\bJ, \bI}.
\end{equation}
If $\alpha = \beta \ne l$, then
\begin{align*}
    X^{\alpha, \beta}_{\bI, \bJ} X^{\alpha, \beta}_{\bJ, \bI}
    &= X^{\alpha, \alpha}_{\bI, \bJ} X^{\alpha,
    \alpha}_{\bJ, \bI} \\
    &= \prod_{s \in \lambda(I^\alpha)} h_{I^\alpha, J^\alpha}(s)
    \epsilon \prod_{s \in \lambda(J^\alpha)} h_{J^\alpha,
    I^\alpha}(s) \epsilon = h_{I^\alpha, J^\alpha}
    h_{J^\alpha, I^\alpha} \epsilon^{|I^\alpha| + |J^\alpha|}.
\end{align*}
By Lemma~\ref{lem:hook=0}, this is equal to zero unless $I^\alpha =
J^\alpha$. We thus need only consider the case where $I^\alpha =
J^\alpha$ for all $\alpha \ne l$.

Note that
\begin{gather*}
    X^{\alpha, \beta}_{\bI, \bJ} = X^{\alpha, \beta}_{\bJ,
    \bJ},\quad X^{\alpha, \beta}_{\bJ, \bI} = X^{\alpha,
    \beta}_{\bI, \bI}, \quad \text{if } \alpha \ne l,\quad
    \text{and} \\
    X^{\alpha, \beta}_{\bI, \bJ} = X^{\alpha, \beta}_{\bI,
    \bI},\quad X^{\alpha, \beta}_{\bJ, \bI} = X^{\alpha,
    \beta}_{\bJ, \bJ}, \quad \text{if } \beta \ne l.
\end{gather*}
Therefore
\[
    \prod_{\begin{array}{cc}
    {\scriptstyle 1 \le \alpha, \beta \le r} \\ {\scriptstyle
    \alpha \ne l \text{ or } \beta \ne l}
    \end{array}} (-1)^{|I^\alpha|} X^{\alpha, \beta}_{\bI, \bJ}
    X^{\alpha, \beta}_{\bJ, \bI} = \prod_{\begin{array}{cc}
    {\scriptstyle 1 \le \alpha, \beta \le r} \\ {\scriptstyle
    \alpha \ne l \text{ or } \beta \ne l}
    \end{array}} (-1)^{|I^\alpha|} X^{\alpha, \beta}_{\bI, \bI}
    X^{\alpha, \beta}_{\bJ, \bJ}.
\]

By Proposition~\ref{prop:ferm-poly-identity-rank-one}, we have
\[
    f_{I^l, J^l} = (-1)^n d_{I^l,J^l} = (-1)^n (-1)^{|I^l|}
    X^{l,l}_{\bI, \bI} X^{l,l}_{\bJ, \bJ},
\]
if $J^l = (-1)^n \psi(k) I^l$ for some $k, n \in \Z$ and $f_{I^l,
J^l}=0$ otherwise.  The result now follows from the definition of
$\psi^l(k)$ and the fact that
\[
    d_{\bI, \bJ} = \prod_{\alpha,\beta = 1}^r (-1)^{|I^\alpha|}
    X^{\alpha, \beta}_{\bI, \bI} X^{\alpha, \beta}_{\bJ, \bJ}.
\]
\end{proof}


\subsection{Proof of the theorem} \label{sec:clifford-proof}

\begin{proof}[Proof of Theorem~\ref{thm:clifford}]
Let $n \in \Z$,
\[
  \bI \in \mathcal{M}_\bc(r,k)^T \quad \text{and} \quad \bJ \in
  \mathcal{M}_{\bc+1_l}(r,k+n-c^l-1)^T.
\]
By Lemma~\ref{lem:structure-constants},
\[
    \langle \Psi^l(n) [\bI], [\bJ] \rangle = \frac{ c_\tnv \left(
    \K_{\bc(\bI), \bc(\bJ)} (r, |\bI|, |\bJ|)_{(\bI, \bJ)}
    \right)}{e_T(\mathcal{T}_\bI^-) e_T(\mathcal{T}_\bJ^+)}
    = \frac{f_{\bI, \bJ}}{d_{\bI, \bJ}}.
\]
The set $\{[\bI]\ |\ \bc(\bI) = \bc,\ |\bI|=n\}$ forms a $\C(b_1,
\dots, b_r, \epsilon)$ basis of the localized equivariant cohomology
of $\M_\bc(r,n)$. Since the above structure coefficients are complex
numbers by Proposition~\ref{prop:ferm-poly-identity-rank-r}, we see
immediately that the operators $\Psi^l(n)$ (and hence their adjoints
$\Psi^l(n)^*$) preserve the subspace $\mathbf{A}$. Let $\varphi$ be
the vector space isomorphism $\mathbf{A} \cong \mathbf{F}$ given by
$\varphi([\bI]) = \bI$.  It follows from
Proposition~\ref{prop:ferm-poly-identity-rank-r} that
$\varphi(\Psi^l(n)[\bI]) = \psi^l(n)\varphi([\bI])$ for all $n \in
\Z$ and $l \in \{1, \dots, r\}$. Since $\Psi^l(n)^*$ and
$\psi^l(n)^*$ are adjoint to $\Psi^l(n)$ and $\psi^l(n)$
respectively, we see that $\varphi(\Psi^l(n)^*[\bI]) =
\psi^l(n)^*\varphi([\bI])$ for all $n \in \Z$ and $l \in \{1, \dots,
r\}$. The result follows.
\end{proof}


\section{Proof of Theorem~\ref{thm:Heisenberg}}
\label{sec:proof-heisenberg}

\subsection{Notation}

In order to simplify notation in what follows, for $r$-colored
semi-infinite monomials $\bI$ and $\bJ$ with $\bc(\bI) =
\bc(\bJ)=\bc$, let
\[
g^{l}_{\bI,\bJ} = \gamma^l_{\bI,\bJ} \cup c_\tnv \left(\K_{\bc,
\bc}(r,|\bI|,|\bJ|)_{(\bI,\bJ)} \right)
\]
in the $T$-equivariant cohomology of a point.  We define
$g^{l}_{\bI, \bJ}$ to be zero when $\bc(\bI) \ne \bc(\bJ)$.  When
non-zero, $g^{l}_{\bI,\bJ}$ is homogeneous of degree $2r(|\bI| +
|\bJ|)$.

Note that when $r=1$, $\lambda(I)_k \ne \lambda(J)_k$ and
$\lambda(I)_i = \lambda(J)_i$ for all $1 \le i < k$, we have
\begin{align*}
    g_{I, J} &\stackrel{\text{def}}{=} g^1_{I,J} \\
    &= \begin{cases}
    (-1)^{|I|-1} \epsilon^{|I|+|J|} h_{J,I} \prod_{s \in \lambda(I),\ s \ne
    (k, \lambda(I)_k)} h_{I,J}(s) & \text{if } \lambda(I)_k >
    \lambda(J)_k, \\
    (-1)^{|I|}\epsilon^{|I|+|J|}
    h_{I,J} \prod_{s \in \lambda(J),\ s \ne (k,\lambda(J)_k)}
    h_{J,I}(s) & \text{if } \lambda(J)_k > \lambda(I)_k. \end{cases}
\end{align*}
Thus, by Lemma~\ref{lem:hook=0}, $g_{I,J}=0$ if and only if $0$
occurs more than once as a relative hook length for the partitions
$\lambda(I)$ and $\lambda(J)$.

Recall that for $\bI \in \M_\bc(r,n_1)^T$, $\bJ \in
\M_{\mathbf{d}}(r,n_2)^T$, we defined
\[
    d_{\bI,\bJ} = e_T(\mathcal{T}_\bI^-)e_T(\mathcal{T}_\bJ^+)
    \in \C[b_1,\hdots,b_r,\epsilon],
\]
and when $r=1$, we have
\[
    d_{I, J} = (-1)^{|I|} h_{I} h_{J} \epsilon^{|I| + |J|}.
\]
We will write $g_{\lambda, \mu}$ and $d_{\lambda, \mu}$ to denote
$g_{I,J}$ and $d_{I,J}$ (respectively) for some semi-infinite
monomials $I$ and $J$ with $\lambda(I) = \lambda$ and
$\lambda(J)=\mu$.  Since $g_{I,J}$ and $d_{I,J}$ are independent of
the charge of $I$ and $J$ , $g_{\lambda, \mu}$ and $d_{\lambda,
\mu}$ are well-defined.

A \emph{$2 \times 2$ square} is a set of the form
\[
    \{(i,j), (i+1,j), (i,j+1), (i+1,j+1)\},\quad i,j \in \N_+.
\]
For $\lambda, \mu \in \mathcal{P}$, we say $\mu - \lambda$ is a
\emph{border strip} if $\lambda \subseteq \mu$ (that is, $\lambda_i
\le \mu_i$ for all $i$) and the following two conditions hold:
\begin{enumerate}
    \item $\mu-\lambda$ contains no $2 \times 2$ square, and
    \item $\mu-\lambda$ is not the disjoint union of two nonempty
    subsets $\nu^1$ and $\nu^2$ such that
    no box in $\nu^1$ shares an edge with a box in $\nu^2$.
\end{enumerate}
The \emph{width} of a border strip is defined to be the number of
columns it occupies minus one.  That is, the width of a border strip
$\mu - \lambda$ is $|\{i\ |\ \lambda_i > \mu_i\}|-1$.  Note that
this is often referred to as the \emph{height} when Young diagrams
are written in English or French notation.


\subsection{Combinatorics of the rank one bosonic Fock space}

\begin{proposition} \label{prop:bosonic-poly-identity-rank-one}
If $|\lambda|<|\mu|$, then
\[
    g_{\lambda, \mu} = \begin{cases}
    \frac{(-1)^{\width (\mu-\lambda)}}{|\mu-\lambda|}d_{\lambda,\mu}
    & \text{if $\mu-\lambda$ is a border strip},\\
    0 & \text{otherwise}. \end{cases}
\]

\begin{proof}
We first prove that $g_{\lambda, \mu}=0$ unless $\lambda \subseteq
\mu$. Since $|\lambda| < |\mu|$, we must have $\mu_k > \lambda_k$
for some positive integer $k$.  Assume $k$ is the smallest such
positive integer.  If $\lambda \not \subseteq \mu$, we must have
$\lambda_j
> \mu_j$ for some $j$.  Again, we take the smallest such $j$.  Then
$h_{\mu, \lambda}((k, \mu_k)) = 0$ and $h_{\lambda, \mu}((j,
\lambda_j))=0$ and thus $g_{\lambda, \mu}=0$.

Now suppose $\lambda \subseteq \mu$.  Let $k$ be the smallest
positive integer such that $\mu_k
> \lambda_k$ and let $s^* = (k, \mu_k)$.  Suppose $\mu-\lambda$
contains a $2 \times 2$ square. Pick one such $2 \times 2$ square
\[
    \{(i,j), (i+1,j), (i,j+1), (i+1,j+1)\}
\]
in $\mu - \lambda$ with the following properties:
\begin{enumerate}
    \item \label{box-prop1} $(i+1,j+2) \not \in \mu$,
    \item \label{box-prop2} $(i-1,j) \not \in \mu - \lambda$.
\end{enumerate}
If property~\eqref{box-prop1} is not satisfied, then
\[
    \{(i,j+1), (i+1,j+1), (i,j+2), (i+1,j+2)\}
\]
is a $2 \times 2$ square in $\mu-\lambda$ and if
property~\eqref{box-prop2} is not satisfied, then
\[
    \{(i-1,j), (i,j), (i-1,j+1), (i,j+1)\}
\]
is a $2 \times 2$ square in $\mu-\lambda$.  Thus, we can always
choose a $2 \times 2$ square in $\mu-\lambda$ satisfying
properties~\eqref{box-prop1} and \eqref{box-prop2} by moving up and
left as necessary.  For such a square, $(i+1,j) \ne s^*$ but
\[
    h_{\mu,\lambda}((i+1,j)) = a_\mu((i+1,j)) + l_\lambda((i+1,j)) + 1 = 1 + (-2) + 1 = 0,
\]
and so $g_{\lambda,\mu}=0$.

Now suppose $\mu - \lambda$ contains no $2 \times 2$ square. It is a
union of subsets $\nu^1, \dots, \nu^p$ such that $\nu^i$ and $\nu^j$
have no edges in common for $i \ne j$. By relabeling if necessary,
we may assume that all boxes in $\nu^i$ are to the left of all boxes
of $\nu^j$ for $i < j$.  By definition $s^*$ is the top left box of
$\nu^1$. If $p \ge 2$, let $(i,j)$ be the top left box of $\nu^2$.
Then $(i-1,j), (i,j+1) \not \in \nu^2$. We also have $(i-1,j),
(i,j+1) \not \in \nu^i$ for $i \ne 2$ since $\nu^2$ shares no edges
with $\nu^i$ for $i \ne 2$.  Thus
\[
    h_{\mu,\lambda}((i,j)) = a_\mu((i,j)) + l_\lambda((i,j)) + 1 = 0 + (-1) + 1 = 0,
\]
and so $g_{\lambda,\mu}=0$.

It remains to consider the case when $\mu-\lambda$ is a border
strip.  We divide the boxes of the partition $\mu$ into subsets as
follows.  Recall that $k$ is the smallest positive integer such that
$\mu_k > \lambda_k$.  Let $m$ be the largest positive integer such
that $\mu_m
> \lambda_m$, let $l$ be the smallest positive integer such that
$\mu^t_l > \lambda^t_l$, and let $n$ be the largest positive integer
such that $\mu^t_n > \lambda^t_n$. Then define (see
Figure~\ref{fig:YD-labels})
\begin{align*}
    A &= \{(i,j) \in \lambda\ |\ (i < k \text{ or } i>m) \text{ and }
        (j<l \text{ or } j>n)\},\\
    B &= \{(i,j) \in \lambda\ |\ i< k,\ l \le j \le n\},\\
    C &= \{(i,j) \in \lambda\ |\ k \le i \le m,\ j < l\}, \text{ and},\\
    D &= \{(i,j) \in \mu\ |\ i \ge k,\ j \ge l\}.
\end{align*}
\begin{figure}
    \includegraphics[width=6cm]{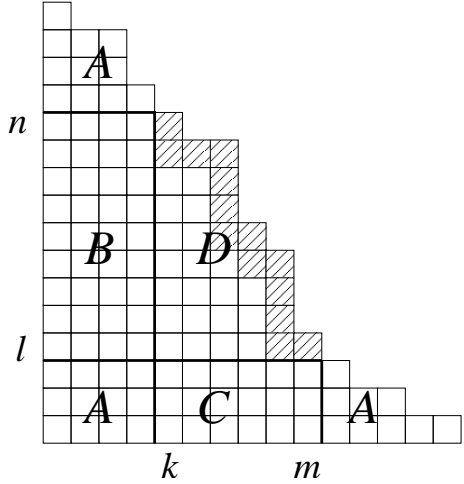}
    \caption{Unshaded boxes correspond to the partition $\lambda$.
    Shaded boxes correspond to the border strip $\mu - \lambda$. \label{fig:YD-labels}}
\end{figure}
We have $\mu = A \sqcup B \sqcup C \sqcup D$ and
\begin{align*}
    h_{\lambda,\mu}(s) &= h_{\mu,\lambda}(s) = h_\lambda(s)
    = h_\mu(s),\quad \text{for } s \in A,\\
    h_{\lambda,\mu}(s) &= h_\mu(s),\quad h_{\mu,\lambda}(s)
    = h_\lambda(s),\quad \text{for } s \in B,\\
    h_{\lambda,\mu}(s) &= h_\lambda(s),\quad h_{\mu,\lambda}(s)
    = h_\mu(s),\quad \text{for } s \in C.
\end{align*}
We now focus our attention on the boxes in $D$.  These boxes form a
partition themselves.  More precisely, if we define
\begin{align*}
    \tilde \lambda_i &= \lambda_{i+k-1}-l+1, \quad 1 \le i \le m+1-k,\\
    \tilde \mu_i &= \mu_{i+k-1}-l+1,\quad 1 \le i \le m+1-k,
\end{align*}
then $\tilde \lambda$ and $\tilde \mu$ are partitions and we have
\begin{align*}
    (i,j) \in \tilde \mu &\iff (i+k-1,j+l-1) \in \mu,\\
    (i,j) \in \tilde \lambda &\iff (i+k-1,j+l-1) \in \lambda,\\
    h_{\tilde \lambda,\tilde \mu}((i,j)) &= h_{\lambda,\mu}((i+k-1,j+l-1)),\\
    h_{\tilde \mu, \tilde \lambda}((i,j)) &= h_{\mu,\lambda}((i+k-1, j+l-1)),\\
    h_{\tilde \lambda}((i,j)) &= h_\lambda((i+k-1, j+l-1)),\\
    h_{\tilde \mu}((i,j)) &= h_\mu((i+k-1, j+l-1)).
\end{align*}
Since $\mu - \lambda$ is a border strip, we have
\begin{align*}
    \lambda_i &= \mu_{i+1}-1,\quad k \le i \le m-1,\\
    \lambda^t_j &= \mu^t_{j+1}-1,\quad l \le j \le n-1.
\end{align*}
Let $\tilde s^* = (1,\tilde \mu_1) = (1, \mu_k - l + 1)$. Let $(i,j)
\in \tilde \lambda$.  Since $\tilde \mu^t_j > \tilde \lambda^t_j$,
we have $(i+1,j) \in \tilde \mu$.  Also
\begin{align*}
    a_{\tilde \mu}((i+1,j)) &= a_{\tilde \lambda}((i,j))+1,\\
    l_{\tilde \lambda}((i+1,j)) &= l_{\tilde \lambda}((i,j)) -1,\\
    l_{\tilde \mu}((i+1,j)) &= l_{\tilde \mu}((i,j)) - 1.
\end{align*}
Thus
\begin{align*}
    h_{\tilde \mu, \tilde \lambda}((i+1,j)) &= a_{\tilde \mu}((i+1,j))
    + l_{\tilde \lambda}((i+1,j)) + 1 \\
    &= a_{\tilde \lambda}((i,j))
    + l_{\tilde \lambda}((i,j)) + 1 = h_{\tilde \lambda}((i,j)),
\end{align*}
and
\begin{align*}
    h_{\tilde \lambda, \tilde \mu}((i,j)) &= a_{\tilde \lambda}((i,j))
    + l_{\tilde \mu}((i,j)) +1 \\
    &= a_{\tilde \mu}((i+1,j)) + l_{\tilde \mu}((i+1,j)) + 1
    = h_{\tilde \mu}((i+1,j)).
\end{align*}
Therefore
\begin{gather*}
    \prod_{(i,j) \in \tilde \lambda} h_{\tilde \mu, \tilde \lambda}((i+1,j))
    = \prod_{(i,j) \in \tilde \lambda} h_{\tilde \lambda}((i,j)) = h_{\tilde \lambda},\\
    \prod_{(i,j) \in \tilde \lambda} h_{\tilde \lambda, \tilde \mu}((i,j))
    = \prod_{(i,j) \in \tilde \lambda} h_{\tilde \mu}((i+1,j)).
\end{gather*}
Now, for $(1,j) \in \tilde \mu$ with $1 \le j < \tilde \mu_1$ (the
second inequality is equivalent to $(1,j) \ne \tilde s^*$), we have
\begin{align*}
    a_{\tilde \mu}((1,j)) &= a_{\tilde \mu}((1,j+1)) + 1,\\
    l_{\tilde \lambda}((1,j)) &= l_{\tilde \mu}((1,j+1)) - 1.
\end{align*}
Thus
\begin{align*}
    h_{\tilde \mu, \tilde \lambda}((1,j)) &= a_{\tilde \mu}((1,j))
    + l_{\tilde \lambda}((1,j)) + 1 \\
    &= a_{\tilde \mu}((1,j+1)) + l_{\tilde \mu}((1,j+1)) + 1 = h_{\tilde \mu}((1,j+1)),
\end{align*}
and so
\[
    \prod_{j\, :\, 1 \le j < \tilde \mu_1} h_{\tilde \mu, \tilde
    \lambda}((1,j)) = \prod_{j\, :\,\ 2 \le j \le \tilde \mu_1} h_{\tilde \mu}((1,j)).
\]
For $1 \le j \le \tilde \mu_1$,
\begin{multline*}
    \prod_{i\, :\, \tilde \lambda^t_i +1 < i \le \tilde \mu^t_i}
    h_{\tilde \mu, \tilde \lambda}((i,j)) = (-1)(-2) \cdots
    (-(\tilde \mu^t_j - \tilde \lambda^t_j - 1)) \\
    = (-1)^{\tilde \mu^t_j
    - \tilde \lambda^t_j -1} (\tilde \mu^t_j - \tilde \lambda^t_j-1)!
\end{multline*}
and
\[
    \prod_{i\, :\, \tilde \lambda^t_i +1 < i \le \tilde \mu^t_i}
    h_{\tilde \mu}((i,j)) = (\tilde \mu^t_j - \tilde \lambda^t_j-1)
    (\tilde \mu^t_j - \tilde \lambda^t_j-2) \cdots (2)(1) =
    (\tilde \mu^t_j - \tilde \lambda^t_j-1)!.
\]
Note that
\[
    \prod_{j=1}^{\tilde \mu_1} (-1)^{\tilde \mu^t_j - \tilde
    \lambda^t_j -1} = (-1)^{\tilde \mu_{\tilde \mu_1}^t -1}
    \prod_{j=1}^{\tilde \mu_1 - 1}
    (-1)^{\tilde \mu^t_j - \tilde \mu^t_{j+1}} = (-1)^{\tilde
    \mu^t_1 - 1} = (-1)^{\width (\mu - \lambda)},
\]
where in the first equality we used the fact that $\tilde
\lambda^t_j = \tilde \mu^t_{j+1} - 1$ for $1 \le j \le \tilde \mu_1
- 1$ and $\tilde \lambda_{\tilde \mu_1}^t=0$. Combining the above
results, we have
\[
    \prod_{s \in \tilde \lambda} h_{\tilde \lambda,
    \tilde \mu}(s) \prod_{s \in \tilde \mu,\ s \ne \tilde s^*}
    h_{\tilde \mu, \tilde \lambda}(s)
    = (-1)^{\width (\mu - \lambda)} \prod_{s \in \tilde \lambda}
    h_{\tilde \lambda}(s) \prod_{s \in \tilde \mu,\ s \ne (1,1)} h_{\tilde
    \mu}(s).
\]
Note that
\[
    h_{\tilde \mu}((1,1)) = |\mu - \lambda|.
\]
Thus
\begin{align*}
    \prod_{s \in D \cap \lambda} h_{\lambda,\mu}(s)
    \prod_{s \in D, s \ne s^*} h_{\mu, \lambda}(s) &= \prod_{s \in \tilde \lambda}
    h_{\tilde \lambda, \tilde \mu}(s) \prod_{s \in \tilde \mu,\ s \ne \tilde s^*}
    h_{\tilde \mu, \tilde \lambda}(s)\\
    &= (-1)^{\width (\mu - \lambda)} \prod_{s \in \tilde \lambda}
    h_{\tilde \lambda}(s) \prod_{s \in \tilde \mu,\ s \ne (1,1)} h_{\tilde
    \mu}(s) \\
    &= (-1)^{\width (\mu - \lambda)} \frac{1}{|\mu-\lambda|}
    \prod_{s \in D \cap \lambda} h_\lambda(s) \prod_{s \in D} h_\mu(s).
\end{align*}
Therefore,
\begin{align*}
    g_{\lambda,\mu} &= (-1)^{|\lambda|}\epsilon^{|\lambda|+|\mu|}
    \prod_{s \in A \cup B \cup C} h_{\lambda, \mu}(s) h_{\mu, \lambda}(s)
    \prod_{s \in D \cap \lambda} h_{\lambda, \mu}(s)
    \prod_{s \in D,\ s \ne s^*} h_{\mu, \lambda}(s) \\
    &= (-1)^{|\lambda|}
    \frac{(-1)^{\width (\mu - \lambda)}}{|\mu-\lambda|} \epsilon^{|\lambda|+|\mu|} \prod_{s \in A \cup B \cup C} h_\lambda(s) h_\mu(s)
    \prod_{s \in D \cap \lambda}
    h_\lambda(s) \prod_{s \in D} h_\mu(s)\\
    &=(-1)^{|\lambda|} \frac{(-1)^{\width (\mu - \lambda)}}{|\mu-\lambda|}
    h_\lambda h_\mu \epsilon^{|\lambda|+|\mu|} \\
    &=\frac{(-1)^{\width (\mu - \lambda)}}{|\mu-\lambda|} d_{\lambda,\mu}
\end{align*}
as desired.
\end{proof}
\end{proposition}

Let $p(-n)$, $n \in \N_+$, be the operator on symmetric functions
given by multiplication by $\frac{1}{n}p_n$, where $p_n$ is the
power-sum symmetric function.

\begin{lemma} \label{lem:schur-function-coefficients}
With respect to the basis of Schur functions, the structure
constants of the operator $p(-n)$ are given by
\[
    \langle p(-n) s_\lambda, s_\mu \rangle = \begin{cases}
    \frac{g_{\lambda,\mu}}{d_{\lambda,\mu}} & \text{if } |\mu -
    \lambda| = n, \\
    0 & \text{otherwise}. \end{cases}
\]
\begin{proof}
By \cite[Example~I.3.11]{Mac95},
\begin{align*}
    p_n s_\lambda = \sum_\mu (-1)^{\width(\mu-\lambda)}
    s_\mu
\end{align*}
where the sum is over all partitions $\mu \supset \lambda$ such that
$\mu - \lambda$ is a border strip of size $n$. Therefore, since the
Schur functions form an orthonormal basis for the symmetric
functions,
\begin{align*}
    \left\langle p(-n) s_\lambda, s_\mu \right\rangle &= \left\langle \frac{1}{n} p_n
    s_\lambda, s_\mu \right\rangle \\
    &= \begin{cases} \frac{1}{n} (-1)^{\width (\mu - \lambda)} &
    \text{if $\mu- \lambda$ is a border strip of size $n$}, \\
    0 & \text{otherwise} \end{cases} \\
    &= \begin{cases}
    \frac{g_{\lambda,\mu}}{d_{\lambda,\mu}} & \text{if } |\mu| -
    |\lambda| = n, \\
    0 & \text{otherwise}. \end{cases}
\end{align*}
\end{proof}
\end{lemma}


\pagebreak

\subsection{Combinatorics of the $r$-colored bosonic Fock space}
\label{sec:bosonic-combinatorics-rank-r}

\begin{proposition} \label{prop:bosonic-poly-identity-rank-r}
If $\bI$ and $\bJ$ are semi-infinite monomials with $\bc(\bI) =
\bc(\bJ)$ and $|\lambda(I^l)| < |\lambda(J^l)|$, then
\[
    g^l_{\bI, \bJ} =
    \frac{(-1)^{\width \left(\lambda\left(J^l\right) -
    \lambda\left(I^l\right)\right)}}{\left|\lambda\left(J^l\right) -
    \lambda\left(I^l\right)\right|} d_{\bI, \bJ}
\]
if $\lambda(I^\alpha) = \lambda(J^\alpha)$ for all $\alpha \ne l$
and $\lambda(J^l) - \lambda(I^l)$ is a border strip. Otherwise
$g^l_{\bI, \bJ}=0$.
\end{proposition}

\begin{proof}
Note that $g^l_{\bI, \bJ}=0$ unless $|I^\alpha|=|J^\alpha|$ for
$\alpha < l$ and $|I^l| \ne |J^l|$ (due to the presence of the
factor $\gamma^l_{\bI, \bJ}$) in which case
\begin{equation} \label{eq:bosonic-polynomial}
    g^l_{\bI, \bJ} = g_{I^l, J^l} \prod_{\begin{array}{cc}
    {\scriptstyle 1 \le \alpha, \beta \le r} \\ {\scriptstyle
    \alpha \ne l \text{ or } \beta \ne l}
    \end{array}} (-1)^{|I^\alpha|} X^{\alpha, \beta}_{\bI, \bJ}
    X^{\alpha, \beta}_{\bJ, \bI},
\end{equation}
where $X^{\alpha, \beta}_{\bI, \bJ}$ is defined in
Section~\ref{sec:fermionic-combinatorics-rank-r}.  As in the proof
of Proposition~\ref{prop:ferm-poly-identity-rank-r}, for all $\alpha
\ne l$ we have $X^{\alpha, \alpha}_{\bI,\bJ}=0$ unless $I^\alpha =
J^\alpha$.  If $I^\alpha=J^\alpha$ for all $\alpha \ne l$, then
\[
    g^l_{\bI, \bJ} = g_{I^l, J^l} \prod_{\begin{array}{cc}
    {\scriptstyle 1 \le \alpha, \beta \le r} \\ {\scriptstyle
    \alpha \ne l \text{ or } \beta \ne l}
    \end{array}} (-1)^{|I^\alpha|} X^{\alpha, \beta}_{\bI, \bI}
    \prod_{\begin{array}{cc}
    {\scriptstyle 1 \le \alpha, \beta \le r} \\ {\scriptstyle
    \alpha \ne l \text{ or } \beta \ne l}
    \end{array}}
    X^{\alpha, \beta}_{\bJ, \bJ}.
\]
By Proposition~\ref{prop:bosonic-poly-identity-rank-one},
\[
    g_{I^l, J^l} = \frac{(-1)^{\width (\lambda(J^l) -
    \lambda(I^l))}}{|\lambda(J^l) - \lambda(I^l)|} d_{I^l,J^l}
    = \frac{(-1)^{\width (\lambda(J^l) - \lambda(I^l))}}
    {|\lambda(J^l) - \lambda(I^l)|} (-1)^{|I^l|} X^{l,l}_{\bI, \bI}
    X^{l,l}_{\bJ, \bJ}
\]
if $\lambda(J^l) - \lambda(I^l)$ is a border strip and is equal to
zero otherwise.  Thus, when $\lambda(J^l) - \lambda(I^l)$ is a
border strip, we have
\[
    g^l_{\bI, \bJ} = \frac{(-1)^{\width (\lambda(J^l) - \lambda(I^l))}}
    {|\lambda(J^l) - \lambda(I^l)|} \prod_{\alpha,
    \beta =1}^r (-1)^{|I^\alpha|} X^{\alpha, \beta}_{\bI, \bI}
    X^{\alpha, \beta}_{\bJ, \bJ}.
\]
The result follows.
\end{proof}

\begin{corollary} \label{cor:schur-function-coefficients-rank-r}
For $n \in \N_+$, we have
\begin{multline} \label{eq:schur-function-coefficients-rank-r}
    \left< p^l(-n) \left(q^{c(I^1)} s_{\lambda(I^1)}, \dots, q^{c(I^r)}
    s_{\lambda(I^r)}\right), \left(q^{c(J^1)} s_{\lambda(J^1)}, \dots, q^{c(J^r)}
    s_{\lambda(J^r)}\right) \right> = \\
    \begin{cases} \frac{g^l_{\bI,
    \bJ}}{d_{\bI, \bJ}} & \text{if } \bc(\bI) = \bc(\bJ),\
    |\lambda(J^l)| - |\lambda(I^l)| = n, \\
    0 & \text{otherwise}. \end{cases}
\end{multline}
\begin{proof}
It follows easily from the definition of the bilinear form and the
operators $p^l(-n)$ that the left side of
\eqref{eq:schur-function-coefficients-rank-r} is equal to zero
unless $\bc(\bI) = \bc(\bJ)$ and $I^\alpha = J^\alpha$ for all
$\alpha \ne l$.  In this case, by
Propositions~\ref{prop:bosonic-poly-identity-rank-one}
and~\ref{prop:bosonic-poly-identity-rank-r},
\[
    \frac{g^l_{\bI, \bJ}}{d_{\bI, \bJ}} = \frac{g_{I^l, J^l}}{d_{I^l,
    J^l}}
\]
and the result follows from
Lemma~\ref{lem:schur-function-coefficients}.
\end{proof}
\end{corollary}


\subsection{Proof of the theorem} \label{sec:heisenberg-proof}

\begin{proof}[Proof of Theorem~\ref{thm:Heisenberg}]
Let $n \in \N_+$, $\bI \in \mathcal{M}_\bc(r,k)^T$ and $\bJ \in
\mathcal{M}_\bc(r,k+n)^T$. By Lemma~\ref{lem:structure-constants},
\[
    \langle P^l(-n) [\bI], [\bJ] \rangle = \frac{\gamma^l_{\bI, \bJ}
    \cup c_\tnv \left( \K_{\bc,\bc}(r, |\bI|, |\bJ|)_{(\bI, \bJ)}
    \right)}{e_T(\mathcal{T}_\bI^-) e_T(\mathcal{T}_\bJ^+)}
    = \frac{g^l_{\bI, \bJ}}{d_{\bI, \bJ}}.
\]
As in Section~\ref{sec:clifford-proof}, we see that the operators
$P^l(-n)$ (and hence their adjoints $P^l(n)$) preserve the subspace
$\mathbf{A}$.  Let $\varphi$ be the vector space isomorphism
$\mathbf{A} \cong \mathbf{B}$ given by $\varphi([\bI]) = (q^{c^1}
s_{\lambda(I^1)}, \dots, q^{c^r} s_{\lambda(I^r)})$.  It follows
from Corollary~\ref{cor:schur-function-coefficients-rank-r} that
$\varphi(P^l(-n)[\bI])) = p^l(-n)\varphi([\bI])$. Since $P^l(n)$ and
$p^l(n)$ are adjoint to $P^l(-n)$ and $p^l(-n)$ respectively, we see
that $\varphi(P^l(n)[\bI]) = p^l(n)\varphi([\bI])$ for all $n \in
\Z$ (the case $n=0$ can be seen directly).  The result follows.
\end{proof}


\section{Vertex operators and geometry} \label{sec:future}

\subsection{A new geometric realization of the boson-fermion
correspondence} \label{sec:BFC-geometric-realization}

From Theorems~\ref{thm:clifford} and~\ref{thm:Heisenberg}, we see
that we have defined actions of the $r$-colored Heisenberg and
Clifford algebras on $\mathbf{A} \subset \bigoplus_{\bc, n}
\mathcal{H}^{2rn}_T (\M_\bc(r,n))$.  In so doing, we obtain a
natural geometrically defined isomorphism between bosonic and
fermionic Fock spaces. Under this isomorphism, the semi-infinite
monomial $\bI$ corresponds to the element $(q^{c(I^1)}
s_{\lambda(I^1)}, \dots, q^{c(I^r)} s_{\lambda(I^r)})$ of bosonic
Fock space since both correspond to the element $[\bI] \in
A_\bc(r,|\bI|)$. This matches up precisely with the classical
(algebraic) boson-fermion correspondence.

The complexes defined in this paper also yield geometric analogues
of the vertex operators appearing in the boson-fermion
correspondence as follows.  Let
\[
    \mathbf{H} = \prod_{\bc, \mathbf{d}, n_1, n_2} \H_T^{2r(n_1+n_2)} (\M_\bc(r, n_1)
    \times \M_{\mathbf{d}}(r, n_2)).
\]
We can think of elements of $\mathbf{H}$ as formal (infinite) linear
combinations of elements of $\H_T^{2r(n_1+n_2)} (\M_\bc(r, n_1)
\times \M_{\mathbf{d}}(r, n_2))$ for $\bc, \mathbf{d} \in \Z^r$,
$n_1, n_2 \in \N$. For $l \in \{1,\dots,r\}$, let $N^l$ be the
operator on $\prod_{\bc, \mathbf{d}, n_1, n_2} \H^*_T(\M_\bc(r,n_1)
\times \M_\mathbf{d}(r,n_2))$ that acts on $\H_T^*(\M_\bc(r,n_1)
\times \M_\bc(r,n_2))$ as $\gamma^l (n_2-n_1) \id$ for $n_1 \ne n_2$
and as $c^l \id$ for $n_1 = n_2$.  It acts on $\H_T^*(\M_\bc(r,n_1)
\times \M_{\mathbf{d}}(r,n_2))$, $\bc \ne \mathbf{d}$, as the
identity (its action on these pieces is actually irrelevant).

Now, the complex \eqref{eq:master-complex} can be considered as a
complex over
\[
    \mathcal{M} \stackrel{\text{def}}{=}
    \bigsqcup_{\bc, \mathbf{d}, n_1, n_2} \M_\bc(r,n_1) \times \M_{\mathbf{d}}(r,n_2).
\]
Let $\K$ denote the vector bundle $\ker \tau/\im \sigma$ on
$\mathcal{M}$ and let $\mathcal{B}$, $\mathcal{B}_-$, and
$\mathcal{B}_+$ be its restrictions to
\begin{gather*}
    \bigsqcup_{\bc,n_1,n_2} \M_\bc(r,n_1) \times \M_\bc(r,n_2),\\
    \bigsqcup_{\bc, n_1, n_2\, :\, n_1 < n_2} \M_\bc(r,n_1) \times
    \M_\bc(r,n_2), \text{ and} \\
    \bigsqcup_{\bc, n_1, n_2\, :\,
    n_1 > n_2} \M_\bc(r,n_1) \times\M_\bc(r,n_2)
\end{gather*}
respectively. Then
\begin{gather}
    N^l c_\tnv(\mathcal{B}) = P^l(0) + \sum_{n \in \Z \setminus \{0\}}
    |n|P^l(n), \label{eq:bosonic-geometric-vertex-op-full} \\
    \gamma^l c_\tnv(\mathcal{B}_-) = \sum_{n>0} P^l(-n), \quad
    \text{\and} \quad -\gamma^l c_\tnv(\mathcal{B}_+) = \sum_{n >0}
    P^l(n) \label{eq:bosonic-geometric-vertex-op-half}
\end{gather}
are elements of $\mathbf{H}$ and are geometric versions of the usual
bosonic vertex operators.  Note that while vertex operators usually
involve a formal variable, this can always be recovered by degree
considerations. Also, when comparing to the presentations in
\cite{K,tKvdL91}, one must remember that the variables $x_i$ used in
the bosonic Fock space there correspond to $p_i/i$.

Now, let $\mathcal{F}^l$ and $\mathcal{F}^{l,*}$ denote the
restriction of $\K$ to
\[
    \bigsqcup_{\bc, n_1, n_2} \M_\bc(r,n_1) \times
    \M_{\bc+1_l}(r,n_2) \quad \text{and} \quad \bigsqcup_{\bc, n_1, n_2} \M_\bc(r,n_1) \times
    \M_{\bc-1_l}(r,n_2)
\]
respectively.  Then
\begin{equation} \label{eq:fermion-geometric-vertex-op}
    c_\tnv(\mathcal{F}^l) = \sum_{n \in \Z} \Psi^l(n) \quad \text{and}
    \quad c_\tnv(\mathcal{F}^{l,*}) = \sum_{n \in \Z} \Psi^l(n)^*
\end{equation}
are elements of $\mathbf{H}$ and are geometric versions of the usual
fermionic vertex operators.

We introduce the usual \emph{normal ordering operator}
\[
    : \Psi^l(k) \Psi^l(j)^* : \stackrel{\text{def}}{=} \begin{cases}
    \Psi^l(k) \Psi^l(j)^* & \text{if $j > 0$}, \\
    -\Psi^l(j)^* \Psi^l(k) & \text{if $j \le 0$}. \end{cases}
\]
Also, for $l \in \{1, \dots, r\}$, define an operator $Q_l :
\mathbf{A} \to \mathbf{A}$ by
\[
    Q_l : A_\bc(r,n) \to A_{\bc+1_l}(r,n),\quad Q_l([\bI]) = [\bJ],
\]
where $\bJ$ is the semi-infinite monomial of charge $\bc + 1_l =
\bc(\bI) + 1_l$ with $\blambda(\bJ) = \blambda(\bI)$.

\begin{proposition} \label{prop:geometric-BFC}
We have
\begin{align}
    N^l c_\tnv(\mathcal{B}) &= \ : c_\tnv(\mathcal{F}^l)
    c_\tnv(\mathcal{F}^{l,*}) :, \label{eq:geom-bosonization} \\
    c_\tnv(\mathcal{F}^l) &= Q_l \exp (\gamma^l c_\tnv(\mathcal{B}_-))
    \exp (\gamma^l c_\tnv(\mathcal{B}_+)), \quad \text{and}
    \label{eq:geom-fermionization-1} \\
    c_\tnv(\mathcal{F}^{l,*}) &= Q_l^{-1} \exp (- \gamma^l
    c_\tnv(\mathcal{B}_-)) \exp (-\gamma^l c_\tnv(\mathcal{B}_+))
    \label{eq:geom-fermionization-2}
\end{align}
as formal operators on $\mathbf{A}$.  More precisely, if we
decompose into blocks according to the decomposition $\mathbf{A} =
\bigoplus_{\bc, n} A_\bc(r,n)$, the above equalities signify
equality of blocks $A_\bc(r,n_1) \to A_{\mathbf{d}}(r,n_2)$ for all
$\bc, \mathbf{d} \in \Z^r$ and $n_1, n_2 \in \N$.
\end{proposition}
\begin{proof}
This follows from the (algebraic) boson-fermion correspondence. See,
for example, \cite[Theorem~1.2]{tKvdL91} and note that our
generators $p^l(m)$ of the oscillator algebra are equal to
$\frac{1}{|m|} \alpha_l(m)$ for $1 \le l \le r$, $m \in \Z \setminus
\{0\}$, where the $\alpha_l(m)$ are the generators used in
\cite{tKvdL91}. One could also produce a geometric proof of this
result in the language of the current paper that mimics the
algebraic proof.
\end{proof}

\begin{remark}
We note the operators $N^l c_\tnv(\mathcal{B})$,
$c_\tnv(\mathcal{F}^l)$, and $c_\tnv(\mathcal{F}^{l,*})$ (more
precisely, their homogeneous components) preserve the integral form
$\mathbf{A}_\Z = \Span_\Z\{[\bI]\}$ of the Fock space $\mathbf{A}$.
\end{remark}

We refer to Equation~\eqref{eq:geom-bosonization} as \emph{geometric
bosonization} and to Equations~\eqref{eq:geom-fermionization-1}
and~\eqref{eq:geom-fermionization-2} as \emph{geometric
fermionization}.  Note that the algebraic boson-fermion
correspondence involves the expression $z^{p^l(0)}$ (or
$z^{\alpha_l(0)}$) where $z$ is the formal variable appearing in the
vertex operators. However, such factors are unnecessary in the
geometric formulation. Their role is played by the relative grading
shift inherent in the definitions of our geometric operators.  More
precisely, it arises from the fact that
\[
    P^l(n)|_{A_\bc(r,k)} \in \H^{2r(2k-n)}_T (\M_\bc(r,k) \times \M_\bc(r,k-n))
\]
while
\begin{gather*}
    \Psi^l(n)|_{A_\bc(r,k)} \in H^{2r(2k+n-c^l-1)}_T (\M_\bc(r,k) \times \M_{\bc+1_l}(r,
    k+n-c^l-1)), \quad \text{and} \\
    \Psi^l(n)^*|_{A_\bc(r,k)} \in H^{2r(2k-n+c^l)}_T (\M_\bc(r,k) \times \M_{\bc-1_l}(r,
    k-n+c^l)).
\end{gather*}


\subsection{Additional vertex operator constructions and future directions}

Vertex operators entered the mathematical literature as a method of
providing explicit constructions of integrable modules for affine
Lie algebras.  In fact, the same integrable representation usually
has several different vertex operator realizations, all of which
should acquire some geometric interpretation using moduli spaces of
sheaves on surfaces.   For the basic level one representation of an
affine Lie algebra $\ag$, vertex operator constructions are
parameterized by conjugacy classes of Heisenberg subalgebras
$\widehat{\mathfrak{h}} \subset \ag$, which correspond bijectively
to conjugacy classes in the finite Weyl group.  In the case of the
affine Lie algebra $\widehat{\mathfrak{gl}}_r$, this means that
there is one vertex operator construction of the basic
representation for each conjugacy class in the symmetric group
$S_r$, and all of these vertex operators have been worked out
explicitly in \cite{tKvdL91} (see also \cite{KP85,Lep85}). For the
conjugacy class of the identity element, the corresponding vertex
operator construction of the basic representation is known as the
homogeneous realization \cite{FK80,Seg81}. This construction
acquires a geometric interpretation in the context of the current
paper using an embedding of $\widehat{\mathfrak{gl}}_r$ into a
completion of the $r$-colored Clifford algebra.  We expect that all
of the vertex operators of the homogeneous realization have a
natural geometric interpretation using the moduli spaces
$\M_\bc(r,n)$ and complexes involving tautological bundles. Indeed,
the $r$-colored Heisenberg algebra action constructed in this paper
defines the action of the homogeneous Heisenberg subalgebra
$\widehat{\mathfrak{h}}_{\mathrm{hom}} \subset
\widehat{\mathfrak{gl}}_r$.

It would be interesting to give geometric realizations of other
vertex operator constructions of the basic representation of
$\widehat{\mathfrak{gl}}_r$, using cyclic group fixed point
components of moduli spaces of framed sheaves on $\mathbb{P}^2$. For
example, the vertex operators of the principal realization
corresponding to the conjugacy class of the Coxeter element should
acquire a geometric interpretation using $\Z_r$-fixed point
components of the moduli space of framed rank 1 sheaves on
$\mathbb{P}^2$, which are Nakajima quiver varieties of type
$\widehat{A}_{r-1}$.  In this context, Nakajima's original Hecke
operators will be identified with homogeneous components of vertex
operators in the principal realization. Geometric constructions of
other realizations of the basic representation would yield, as a
corollary, a geometric interpretation of the vertex operators needed
for level-rank duality.  We hope to say more about this in a future
paper.


\bibliographystyle{abbrv}
\bibliography{biblist}

\end{document}